\documentclass[10pt]{amsart}
\usepackage{amsmath,amssymb,amsfonts,amsthm,amsopn}
\usepackage{latexsym}
\usepackage{mathrsfs}
\usepackage{geometry}
\usepackage{xcolor}
\usepackage{color, colortbl}
\usepackage{pb-diagram}
\usepackage[title]{appendix}
\usepackage{tikz-cd}
\usepackage{tikz}
\usepackage{appendix}
\usepackage{mathtools}
\usepackage{booktabs}
\usepackage{enumerate}
\usepackage{hyperref}
\usepackage{multirow}
\usepackage{comment}
\usepackage{xcolor}
\usetikzlibrary{arrows.meta}
\usetikzlibrary{angles, quotes, arrows.meta, patterns, decorations.pathreplacing, calc, intersections}

\newcommand{\Imab}{\Im(a/b)}

\mathtoolsset{showonlyrefs}

\usepackage{multirow}

\newtheorem{theorem}{Theorem}[section]
\newtheorem{lemma}[theorem]{Lemma}
\newtheorem{corollary}[theorem]{Corollary}
\newtheorem{proposition}[theorem]{Proposition}
\newtheorem{definition}[theorem]{Definition}

\newtheorem{remark}[theorem]{Remark}
\theoremstyle{definition}
\newtheorem{example}[theorem]{Example}

\newcommand{\beqa}{\begin{eqnarray*}}
	\newcommand{\eeqa}{\end{eqnarray*}}

\DeclareMathOperator*{\Sp}{Sp}
\DeclareMathOperator*{\Mp}{Mp}
\DeclareMathOperator*{\Sym}{Sym}

\newcommand{\field}[1]{\mathbb{#1}}
\newcommand{\bR}{\field{R}}        
\newcommand{\bZ}{\field{Z}}        
\newcommand{\bC}{\field{C}}        
\newcommand{\bF}{\field{F}} 
\newcommand{\bT}{\field{T}}        %
\newcommand{\bH}{\field{H}}        %

\def\la{\lambda}

\def\Fur{\mathscr{F}}              
\def\cS{\mathcal{S}}
\def\cD{\mathcal{D}}

\def\cQ{\mathcal{Q}}

\def\cC{\mathcal{C}}

\def\rd{\bR^d}
\def\cd{\bC^d}

\def\R{\right)}

\def\<{\left<}
\def\>{\right>}

\def\ud{\, \mathrm{d}}

\def\mv1{M_v^1}

\def\mn{(m,n)}
\def\mn'{(m',n')}

\newcommand{\norm}[1]{\lVert#1\rVert}

\def\c{\hfill\break}

\def\w{\mathrm{w}}

\def\R{\mathbb{R}}
\def\Ren{\mathbb{R}^d}

\def\Fur{\mathcal{F}}

\def\Sn2{S_{2}(L^{2}(\Ren))}
\def\S1{S_{1}(L^{2}(\Ren))}
\def\sig00{\sigma_{0,0}}

\def\la{\langle}
\def\ra{\rangle}

\begin{document}
	
	\begin{abstract}
         Uncertainty principles are fundamental constraints expressing the impossibility of a signal being simultaneously concentrated in both time and frequency. 
         The classical Hardy uncertainty principle makes this precise in terms of Gaussian decay: if the joint Gaussian decay of a function and its Fourier transform exceeds a certain threshold, the function must be identically zero. 
         In this work we extend this principle to the metaplectic (or oscillator) semigroup, which contains the Fourier transform and the propagators of dissipative quadratic evolutions.
         In the dissipative setting, the classical extremal picture breaks down: the critical threshold may fail to admit any non-zero extremizer. When a non-zero critical extremizer does exist, it is a bounded quadratic oscillation rather than an integrable Gaussian.
         We give a complete and sharp description of all regimes, and derive dynamic versions of Hardy's uncertainty principle for concrete Schrödinger problems.
	\end{abstract}
	
	\title[Hardy's uncertainty principle for the oscillator semigroup]{Hardy's uncertainty principle for the oscillator semigroup, with applications to quadratic evolutions}

    \author{Elena Cordero}
	\address{Universit\`a di Torino, Dipartimento di Matematica}
	\email{elena.cordero@unito.it}
    
	\author{Gianluca Giacchi}
	\address{Universit\`a della Svizzera Italiana}
	\email{gianluca.giacchi@usi.ch}

	\author{Eugenia Malinnikova}
    \address{Stanford University $\&$ \newline
    Norwegian University of Science and Technology}
    \email{eugeniam@stanford.edu}
	
	\thanks{
    Elena Cordero and Gianluca Giacchi have been supported by the Gruppo Nazionale per l’Analisi Matematica, la Probabilità e le loro Applicazioni (GNAMPA) of the Istituto Nazionale di Alta Matematica (INdAM).
    Gianluca Giacchi has been supported by the SNSF starting grant ``Multiresolution methods for unstructured data" (TMSGI2 211684). Eugenia Malinnikova was supported by NSF grant DMS-2247185.
    }

\maketitle

\section{Introduction}
\subsection{Classical Hardy's uncertainty principle}
Uncertainty principles are the neatest expressions of the dualism between localization in time and frequency: a non-zero signal cannot be simultaneously well-localized along with its Fourier transform, and different measures of localization, or concentration, yield different uncertainty principles. In this work, we are mainly concerned with Hardy-type uncertainty principles, where the notion of concentration is given in terms of Gaussian decay: a non-zero function and its Fourier transform cannot decay both arbitrarily fast at infinity. The very first formulation of this phenomenon was given by Hardy in 1933 for functions in $L^2(\bR)$, cf. \cite{Hardy}.
    \begin{theorem}[Classical Hardy’s uncertainty principle]\label{intro.HardyClassical}
	Let $f \in L^2(\mathbb{R})$ and let $\delta,\gamma > 0$. Assume that
	\begin{equation}\label{eq:Hardy_decay_f}
		|f(x)| \lesssim e^{-\pi \delta x^2}, \qquad 
        |\widehat{f}(\xi)| \lesssim e^{-\pi \gamma \xi^2}, 
		\qquad x,\xi \in \mathbb{R}.
	\end{equation}
	\begin{itemize}
		\item If $\gamma\delta>1$, then $f \equiv 0$.
		\item If $\gamma\delta=1$, then
		\begin{equation}\label{eq:Hardy_extremizer}
			f(x) = {C}  e^{-\pi \delta x^2}
		\end{equation}
		for some constant ${C} \in \mathbb{C}$.
        \item If $\gamma\delta<1$, the Gaussian estimates in \eqref{eq:Hardy_decay_f} are satisfied by a dense subspace of $L^2(\bR)$.
	\end{itemize}
\end{theorem}
Notably, Gaussian states are the optimizers in Theorem \ref{intro.HardyClassical}, a recurring feature when it comes to uncertainty principles.
The original result by Hardy was later extended in several directions. Sitaram, Sundari and Thangavelu proved it for functions in several variables \cite{SST} by means of the Radon transform, while Morgan replaced Gaussians with functions in the form $e^{-\pi \alpha|x|^p}$, \cite{Morgan}, thereby allowing more general exponential decays. Distributional versions of Hardy's uncertainty principle can be found in the survey by Bonami and Demange \cite{BDsurvey}. Dynamic versions of Hardy's uncertainty principle for Schr\"odinger evolutions, including equations with potentials, were developed in the notable contributions by Escauriaza, Kenig, Ponce and Vega \cite{EKPV2008,EKPV2010}. Gaussian decay for harmonic oscillators and related models was investigated by Cassano and Fanelli in \cite{CassanoFanelli1}.
A geometric perspective on quantum indeterminacy and its connections with uncertainty principles is developed in the recent preprint by de Gosson \cite{de2026quantum}. 
Finally, D'Ancona and Fiorletta \cite{Diego1} establish dynamical Amrein-Berthier uncertainty principles for magnetic Schr\"odinger flows, subsequently extended by D'Ancona, Fiorletta and Nicola \cite{Diego2} through a geometric nonrefocusing criterion.

\subsection{Hardy's uncertainty principle for \texorpdfstring{$\Mp(\bR)$}{}}
In the present contribution, we are interested in extending classical Hardy's uncertainty principle to a class of operators containing propagators of quadratic evolution equations as particular instances, following the work \cite{FBM2021}. There, the following dynamic Hardy-type uncertainty principle for the heat equation is discussed:
\begin{theorem}\label{intro.ThmMalinnikova}
    Let $u\in\mathcal{C}^1([0,t_0],H^2(\bR))$ be a solution of the heat problem
    \begin{equation}\label{heat}
        \partial_t u=\partial_x^2 u \quad \text{on $[0,t_0]\times\bR$}, \qquad u(0,\cdot)=u_0\in L^1(\bR).
    \end{equation}
    If 
    \begin{equation}
        |u(t_0,x)|\lesssim e^{-\pi \delta x^2},
    \end{equation}
    for some $\delta>0$, and $\delta\geq 1/(4\pi t_0)$, then $u=0$.
\end{theorem}
Observe that, differently from Theorem \ref{intro.HardyClassical}, the critical exponent $\delta= 1/(4\pi t_0)$ also yields $u=0$.
In general, the propagator $e^{-2\pi ita^\w(x,D)}$ of the Schr\"odinger problem
\begin{equation}\label{intro.Schrodinger}
    \begin{cases}
        \frac{i}{2\pi}\partial_tu=a^{\w}(x,D)u, & \text{$t>0$, $x\in\bR$},\\
        u(0,\cdot)=u_0\in\cS(\bR),
    \end{cases}
\end{equation}
where $a^\w(x,D)$ is the Weyl quantization of a complex quadratic form $a$ with negative semi-definite imaginary part, as detailed in Section \ref{sec:Examples}, defines a one-parameter semigroup of contractions. In turn, these operators are particular instances of so-called {\em metaplectic operators}. 
Let us briefly discuss their definition. The group of $2\times2$ complex symplectic matrices \begin{equation}\label{intro.blockU}
    S=\begin{pmatrix}
        a & b\\
        c & d
    \end{pmatrix}, \qquad ad-bc=1,
\end{equation}
admits a subsemigroup $\Sp_+(\bC)$ of matrices such that the associated fractional linear map
\begin{equation}
    \mathcal{C}(S)(\lambda)=(c+d\lambda)(a+b\lambda)^{-1},
\end{equation}
is well-defined and holomorphic on the open upper-half plane $\mathbb U^\circ$, and
\begin{equation}
    \lambda\in \mathbb{U}^\circ \qquad \Longrightarrow\qquad \mathcal{C}(S)(\lambda)\in\mathbb{U}^\circ,
\end{equation}
see, e.g., \cite{Folland}.
The {\em metaplectic semigroup} (or, {\em oscillator semigroup}) $\Mp_+(\bC)$ is the two-fold cover of $\Sp_+(\bC)$, see Section \ref{sec:prerequisites} below for the details. 
The covering $\widehat S\in\Mp_+(\bC)\mapsto S\in\Sp_+(\bC)$ is also known as {\em projection} and it is a homomorphism.
The group of units $\Mp_+(\bC)^\times$ of $\Mp_+(\bC)$ is the {\em metaplectic group} $\Mp(\bR)$. 
This reflects the fact that if $a^\w(x,D)$ in \eqref{intro.Schrodinger} is the quantization of a real quadratic form, the corresponding propagator $e^{-2\pi ita^\w(x,D)}$ defines a one-parameter subgroup of metaplectic operators, and \eqref{intro.Schrodinger} is indeed well-posed for every $t\in\bR$. 

To the best of our knowledge, Hardy's uncertainty principle for operators in
$\Mp(\bR)$ was first established by Zhang \cite{Hardy1}. 
The unitary case is by now well understood. Costa Dias, de Gosson and Prata
proved Hardy's uncertainty principle for the so-called \emph{quadratic Fourier
transforms} in dimension $d\geq1$, see \cite[Theorem 24]{DDP2024}. Combining
their techniques with the integral representation formulas of ter Morsche and
Oonincx \cite{TMO}, the result was extended to the whole metaplectic group in
\cite{CGM2025}, where the multivariate picture is shown to be considerably
more varied, cf. \cite[Theorem 4.2]{CGM2025}.
On the other hand, metaplectic operators with $b=0$ do not satisfy Hardy-type
uncertainty principles. Hence, in the one-dimensional setting, quadratic
Fourier transforms exhaust the class of metaplectic operators for which such
uncertainty principles hold, and Theorem \ref{intro.thm13} provides a complete
characterization.

We recall the one-dimensional
statement, which will be the natural term of comparison for our results.

\begin{theorem}[Hardy's uncertainty principle for $\Mp(\bR)$]\label{intro.thm13}
    Let $\widehat U\in\Mp(\bR)$ with projection $U\in\Sp(\bR)$ with real entries as in \eqref{intro.blockU}.
    Assume $b\neq0$ and
    \begin{align}\label{intro.decayHardy}
        |f(x)|\lesssim e^{-\pi \delta x^2},\qquad |\widehat Uf(x)|\lesssim e^{-\pi\gamma x^2}, \qquad x\in\bR,
    \end{align}
    for some $\gamma,\delta>0$.
    \begin{itemize}
        \item If $\delta\gamma>b^{-2}$, then $f=0$.
        \item If $\delta\gamma=b^{-2}$, then $f(x)={C} e^{-\pi (\delta+i a/b) x^2}$, for some $C\in\bC$.
        \item If $\delta\gamma<b^{-2}$, there exists a dense subset of functions in $L^2(\bR)$ satisfying \eqref{intro.decayHardy}.
    \end{itemize}
\end{theorem}

Finally, time-frequency methods are adopted in \cite{BDJ203,CGPT26} to relate decay of
functions or distributions to decay estimates for the Wigner transform. In
\cite{grochenig2025more} Hardy and Beurling type uncertainty principles for metaplectic Wigner distributions in any dimension are derived from results on partial short-time Fourier transforms. The authors also obtain a supercritical Hardy type theorem with scalar Gaussian majorant.
For the directional Hardy theorem for unitary metaplectic operators,
we refer to \cite{CGM2025},
and for further anisotropic uncertainty principles for unitary metaplectic operators in any dimension  to \cite{CGP2026}.

The subcritical case of Theorem~\ref{intro.thm13} admits a quantitative description in terms of Hermite coefficients. The following statement follows
from Vemuri's result
\cite[Theorem 2.2 and its proof]{Vemuri2008}
by a change of scale and multiplication by a quadratic phase. For $n=0,1,2,\ldots,$ let
\[
    \psi_n(x)
    =\frac{2^{1/4}}{(2^n n!)^{1/2}}
      H_n(\sqrt{2\pi}\,x)e^{-\pi x^2},
    \qquad {\text{where}}\quad
    H_n(t)=(-1)^n e^{t^2}
           \frac{\ud^n}{\ud t^n}e^{-t^2}
\]
is the $n$-th Hermite polynomial.

\begin{theorem}[Metaplectic form of Vemuri's estimate]
\label{intro.Vemuri}
    Let $\widehat U\in\Mp(\bR)$ be as in Theorem \ref{intro.thm13} and let $f$ satisfy \eqref{intro.decayHardy}, where $\delta\gamma<b^{-2}$, $b\neq 0$.
     Set
    \[
        \ell=\left(\frac{b^2\gamma}{\delta}\right)^{1/4},
        \qquad
        \psi_n^{U,\delta,\gamma}(x)
        =\ell^{-1/2}e^{-i\pi(a/b)x^2}\psi_n(x/\ell).
    \]
    Then $(\psi_n^{U,\delta,\gamma})_{n\geq0}$ is an
    orthonormal basis of $L^2(\bR)$, and
    \[
        \big|\langle f,\psi_n^{U,\delta,\gamma}\rangle\big|
        \lesssim  (n+1)^{-1/4}
        \left(
            \frac{1-|b|\sqrt{\delta\gamma}}
                 {1+|b|\sqrt{\delta\gamma}}
        \right)^{n/4},
        \qquad n=0,1,2,\ldots.
    \]
\end{theorem}

\subsection{Present contributions}
The metaplectic uncertainty principles recalled above concern the real symplectic group. Passing to the metaplectic semigroup $\Mp_+(\bC)$ is not a technical
extension of the unitary theory: dissipation alters the nature of the problem
in three distinct ways. First, the kernel of $\widehat S$ may already supply a Gaussian decay factor on its output, which must be accounted for in the formulation of a Hardy-type hypothesis, see Proposition~\ref{propFreeDecay} below. Second, the
critical threshold is no longer the algebraic quantity $b^{-2}$, but is
governed by $\Im(a/b)$ and $|b|^{-2}$, that is, by the geometry of the
Hamiltonian flow. Third, and most importantly, the extremal theory is no
longer rigid: saturating the critical threshold no longer guarantees the existence of
an optimizer.
We carry out this analysis in dimension one, where the interplay of these
effects admits a complete and sharp description, and where the answer already
depends on algebraic conditions on the entries of the symplectic matrix that
have no unitary counterpart.

While the metaplectic semigroup traces its origins back to the 1970s, through the early works of Brunet and Kramer \cite{BK1,BK2,Brunet}, Olshanskii \cite{Olshanskii}, Howe \cite{Howe}, and Hörmander \cite{Hormander1,Hormander3}, its structure in relation to the double covering of $\Sp_+(\bC)$ remains far from fully explicit. 
In most treatments, attention is indeed restricted to those operators that arise as propagators of problems \eqref{intro.Schrodinger}.

\subsubsection{Integral representation formulas}
The first contribution of this work is a more operator-theoretic characterization of $\Sp_+(\mathbb{C})$, formulated in terms of the positive semi-definiteness of a certain quadratic form.
\begin{theorem}\label{intro.Main-1}
      Let $S\in\Sp(\bC)$ have entries as in \eqref{intro.blockU}. 
      \begin{enumerate}[(i)]
          \item If $b=0$, then $S\in\Sp_+(\bC)$ if and only if
      \begin{equation}
          a\in\bR^\times \qquad \text{and} \qquad \Im\big(c/a\big)\geq0.
      \end{equation}
      \item If $b\neq0$, then $S\in\Sp_+(\bC)$ if and only if
      \begin{equation}\label{intro.posDefQF}
             \Im\begin{pmatrix}
                  d/b & -1/b\\
                  -1/b & a/b
              \end{pmatrix}\geq0.
          \end{equation}
      \end{enumerate}
  \end{theorem}
Thanks to Theorem \ref{intro.Main-1}, we are able to establish the following representation formulas for operators in $\Mp_+(\bC)$.
\begin{theorem}[Integral representation formulas for $\Mp_+(\bC)$]\label{intro.main-2}
      Let $\widehat S\in\Mp_+(\bC)$ have projection $S\in\Sp_+(\bC)$ with entries as in \eqref{intro.blockU}.
      \begin{enumerate}[(i)]
          \item If $b=0$, then
          \begin{equation}\label{intro.Nonfree}
             \widehat Sf(x)=\zeta |a|^{-1/2}e^{i\pi cx^2/a}f\big(x/a\big), \qquad f\in L^2(\bR).
          \end{equation}
          for some $\zeta\in\bT$ {\em (phase factor)}.
          \item If $b\neq0$, then
          \begin{equation}\label{into.freeMetapOps}
              \widehat Sf(x)=\zeta |b|^{-1/2}e^{i\pi dx^2/b}\int_{\bR}f(y)e^{i\pi ay^2/b}e^{-2\pi ixy/b}\ud y, \qquad f\in \cS(\bR).
          \end{equation}
          for some $\zeta\in\bT$.
      \end{enumerate}
  \end{theorem}
  To simplify the discussion, throughout this work we shall use the phraseology ``up to a phase factor" regarding equalities that are true up to the correct choice of a multiplicative phase factor, whilst leaving this constant implicit.

    Observe that, up to a constant, formula \eqref{into.freeMetapOps} includes all the operators with complex Gaussian kernels
    \begin{equation}
        Tf(x)=\int_{\bR}f(y)e^{2\pi iq(x,y)}\ud y,
    \end{equation}
  where $q(x,y)=\frac 1 2 Py^2-Exy+\frac 1 2 Qx^2$, $P,Q\in\bC$ and $E\in\bC^\times$, is a complex quadratic form with positive semi-definite imaginary part. In this case, $T=C\widehat S$ for some $C\in\bC$ and
  \begin{equation}
      S=\begin{pmatrix}
          PE^{-1} & E^{-1}\\
          QPE^{-1}-E & QE^{-1}
      \end{pmatrix}.
  \end{equation}
  Observe that $S\in\Sp_+(\bC)$ if and only if $\Im(q)\geq0$, in view of \eqref{intro.posDefQF}, i.e., exactly the positive semi-definiteness assumption on $q$.
  
  Differently from generalized Mehler formulas \cite{Hormander3}, the integral representation \eqref{into.freeMetapOps} is given directly in terms of the entries of a symplectic matrix. 
  
   Concerning the convergence of \eqref{into.freeMetapOps}, we remark that $\Im(a/b),\Im(d/b)\geq0$ for $S\in\Sp_+(\bC)$, due to \eqref{intro.posDefQF}, and therefore the integral therein converges whenever $f\in\cS(\bR)$, and it defines a function in $L^2(\bR)$. 
   Furthermore, if $\Im(a/b)>0$, then $f\in L^p(\bR)$, $1\leq p\leq \infty$, is also sufficient for the integral to converge absolutely.
   
   \subsubsection{The main result} Exactly as in the unitary setting, if $b=0$, no uncertainty principle can hold, since \eqref{intro.Nonfree} shows that $\widehat S$ maps compactly supported functions to compactly supported functions. 
    In the case $b\neq0$, the integral representation formula \eqref{into.freeMetapOps} yields the following result.
    Since $\Im(a/b)\geq0$ for every $S\in\Sp_+(\bC)$ with $b\neq0$, by \eqref{intro.posDefQF}, the following two cases are mutually exclusive and exhaustive.
    
    \begin{theorem}[Hardy's uncertainty principle for $\Mp_+(\bC)$]\label{thm-main3}
    Let $\widehat S\in\Mp_+(\bC)$ have projection $S$ as in \eqref{blockS} with $b\neq0$. 
    \begin{enumerate}[(i)]
    \item \emph{Case $\Im(a/b)=0$.}
      Suppose that $f\in L^2(\bR)$ satisfies
      \begin{equation}\label{intro.decaySf2}
          |f(x)| \lesssim e^{-\pi\delta x^2}, \qquad
          |\widehat{S}f(x)| \lesssim e^{-\pi\big(\gamma+\Im(d/b)\big)x^2},
          \qquad x\in\bR,
      \end{equation}
      for some $\gamma>0$ and $\delta>0$.
      \begin{itemize}
          \item If $\gamma\delta > b^{-2}$, then $f = 0$.
          \item If $\gamma\delta = b^{-2}$, then
                \begin{equation}
                    f(x) = {C} e^{-\pi\big(\delta + ia/b\big)\,x^2}, 
                \end{equation} 
                for some ${C}\in\bC$.
          \item If $\gamma\delta < b^{-2}$, a dense subset of $L^2(\bR)$
                satisfies both decay conditions.
      \end{itemize}
    \item \emph{Case $\Im(a/b)>0$.}
    Suppose that $f\in L^p(\bR)$, $1\leq p\leq\infty$, satisfies 
        \begin{equation}\label{intro.decaySf1}
            |\widehat Sf(x)|\lesssim e^{-\pi\big(\gamma+\Im(d/b)-\Im(a/b)^{-1}\Im(b)^{2}/|b|^4\big)x^2}, \qquad x\in\bR,
        \end{equation}
        for some $\gamma\geq 0$.
        \begin{itemize}
            \item If $\gamma\Im(a/b) > |b|^{-2}$, then $f=0$.
            \item If $\gamma\Im(a/b)= |b|^{-2}$, then $f=0$, unless $p=\infty$ and $\Re(b)\neq0$, in which case
            \begin{equation}\label{intro.def-fb20}
                f(x)=C e^{-i\pi \big(\Re(a)/\Re(b)\big)x^2}
            \end{equation}
            for some $C\in\bC$.
      \item If $\gamma \Im(a/b)<|b|^{-2}$, the space of functions satisfying \eqref{intro.decaySf1} is dense in $L^2(\bR)$.
        \end{itemize}
    \end{enumerate}
\end{theorem}
    The inequalities in \eqref{intro.decaySf2} and \eqref{intro.decaySf1} hold pointwise and Item $(i)$ recaptures Theorem \ref{intro.thm13}. If $\Imab=0$, then $b\in\bR$ (see Corollary \ref{corbReal}), and it is evident from \eqref{into.freeMetapOps} that in this case
    \begin{equation}
        \widehat Sf=|b|^{-1/2}e^{i\pi d(\cdot)^2/b}\Fur\big(f e^{i\pi a(\cdot)^2/b}\big)\big(\cdot/b\big)\in \cC^\infty(\bR),
    \end{equation}
    so that the Gaussian decay of $f$ implies $\widehat Sf\in\mathcal{C}^\infty(\bR)$. Conversely, the decay condition for $\widehat Sf$ in \eqref{intro.decaySf2} implies that $\widehat Sf e^{-i\pi d(\cdot)^2/b}$ has Gaussian decay, whence $$f=|b|^{1/2}e^{-i\pi a(\cdot)^2/b}\Fur^{-1}\big(\widehat Sf(b\cdot) e^{-i\pi db(\cdot)^2}\big)\in\cC^\infty(\bR).$$ 
    In Item $(ii)$ of Theorem \ref{thm-main3}, the integral representation defines a continuous representative of $\widehat Sf$, as shown in Proposition \ref{propFreeDecay}. We use this representative for pointwise evaluations.

  The subcritical regime in Theorem \ref{thm-main3} $(ii)$ admits a quantitative description in terms of Hermite coefficients similar to Theorem \ref{intro.Vemuri}. We adapt the argument of Vemuri and obtain an analogous estimate for a normalized output of the metaplectic semigroup.

    \begin{theorem}[Subcritical Hermite coefficient estimates]
\label{intro.subcriticalHermite}
    Let $f$ and $\widehat S$ satisfy the assumptions of
    Theorem \ref{thm-main3} $(ii)$, with
    $0<\gamma<\gamma_0=(\Im (a/b))^{-1}|b|^{-2}$. Set
    \[
       \omega_S=-\frac{i}{b}\bigl(d+\gamma_0\Im b\bigr),
    \qquad
    \ell=(\gamma\gamma_0)^{-1/4},\quad  
        \qquad
        h(x)=\ell^{1/2}
             e^{\pi\omega_S\ell^2x^2}
             (\widehat Sf)(\ell x).
    \]
    Then $h\in L^2(\bR)$ and
    \[
        |\langle h,\psi_n\rangle|
        \lesssim (n+1)^{-1/4}
        \left(
            \frac{1-\sqrt{\gamma/\gamma_0}}
                 {1+\sqrt{\gamma/\gamma_0}}
        \right)^{n/4},
        \qquad n=0,1,2,\ldots.
    \]
\end{theorem}

\begin{remark}
Under the assumptions of Theorem~\ref{intro.subcriticalHermite},
suppose in addition that $b$ is real. Set $r=|b|/\ell$ and
\[
    g(x)=r^{1/2}e^{i\pi(a/b)r^2x^2}f(rx).
\]
Since $\omega_S=-id/b$, the integral formula
\eqref{into.freeMetapOps} gives, up to a constant phase,
\[
    h(x)=\widehat g\bigl(\operatorname{sgn}(b)x\bigr).
\]
The Hermite functions are Fourier eigenfunctions and have
definite parity. Hence
\[
    |\langle g,\psi_n\rangle|
    =|\langle h,\psi_n\rangle|,
\]
so the estimate of Theorem~\ref{intro.subcriticalHermite}
also holds for the Hermite coefficients of $g$.
When $p=\infty$, this follows directly from Vemuri's theorem:
both $g$ and $\widehat g$ satisfy the Gaussian bound with
exponent $\sqrt{\gamma/\gamma_0}$.
\end{remark}

    \subsubsection{The breakdown of the extremal theory}
    The critical exponent case of the classical Hardy uncertainty principle, corresponding to $\gamma \delta=b^{-2}$ in Theorem \ref{intro.thm13}, is clean and complete: the only functions saturating both Gaussian decay bounds are chirped Gaussians  $f(x) = {C} e^{-\pi(\delta+ia/b) x^2}$, unique up to a scalar multiple. In particular, a non-trivial optimizer always exists at the critical threshold. 
    When $\Imab>0$, the assumption $f\in L^p(\bR)$ appears as a natural sufficient condition for the integral in \eqref{into.freeMetapOps} to converge.
    In this case, one of the most striking features of our main results is that the picture discussed above in the critical regime breaks down in the metaplectic semigroup setting, revealing a more nuanced situation.
    At the critical threshold $\gamma\Im(a/b)=|b|^{-2}$, the existence of non-trivial extremizers is no longer guaranteed. As stated in Theorem \ref{thm-main3} $(ii)$, a necessary and sufficient condition for their existence in $L^p(\bR)$ is that $\Re(b)\neq0$ and $p=\infty$. 
    Otherwise, the critical threshold is reached and yet the only admissible Lebesgue function is $f=0$, a phenomenon with no counterpart in the classical theory. 
    
    \subsection{Comments on the main result}
    We now comment on Theorem \ref{thm-main3}.
    \begin{enumerate}
        \item The Gaussian exponents in \eqref{intro.decaySf2} and \eqref{intro.decaySf1} are not imposed for technical necessity. 
    Rather, they account for the natural decay arising intrinsically from the dissipative effect of $\widehat S$, as detailed in Proposition \ref{propFreeDecay} in Section \ref{secHUP}. On top of this, their relation to $\widehat S$ is not as obscure as it may seem, for instance $\Im(d/b)-\Imab^{-1}\Im(b)^2/|b|^4$ is, up to a positive constant, the  determinant of the imaginary part of the quadratic form defining the kernel of $\widehat S$, see \eqref{intro.posDefQF} above.
    \item The case $\Im(a/b)=0$ closely mirrors Theorem \ref{intro.thm13}, which in fact appears as a particular case of Theorem \ref{thm-main3} $(i)$. When $\Imab=0$, the factor $e^{i\pi a x^2/b}$ is a chirp rather than a Gaussian. Consequently, the Gaussian decay assumption on $f$ in \eqref{intro.decaySf2} becomes essential in the proof of Theorem \ref{thm-main3} to counterbalance the lack of the Gaussian factor in \eqref{into.freeMetapOps}.
    \item  An analogue of Theorem \ref{thm-main3} can be obtained for metaplectic operators $\widehat S$ having projection $S$ with $a\neq0$ (and $b\neq0$, to avoid triviality). In this case, the integral representation of $\widehat S$ yields Hardy's uncertainty principle for distributions in the spaces $\Fur L^p(\bR)$, see Theorem \ref{thm-main4} below. 
    \item The critical phase in \eqref{intro.def-fb20} has a direct interpretation in terms of the polar decomposition $\widehat S=\widehat U\widehat Z$, where $\widehat U\in\Mp(\bR)$ is unitary and $\widehat Z$ is positive, see \cite[Theorem 4.4 and Corollary 4.5(ii)]{GRT2026}.
    In dimension one, the projection of the positive factor has real part equal to a positive scalar multiple of the identity. Consequently,
    \[
        U=\begin{pmatrix}
            \alpha&\beta\\
            \omega&\rho
        \end{pmatrix}
        =\frac{\Re(S)}{\sqrt{\det(\Re(S))}}.
    \]
    Thus $\Re(b)\neq0$ if and only if $\beta\neq0$, and in this case
    \[
        \frac{\Re(a)}{\Re(b)}=\frac{\alpha}{\beta}.
    \]
    In the regime $\Im(a/b)>0$, non-zero critical extremizers in $L^p(\bR)$ therefore exist precisely when $p=\infty$ and the unitary polar factor $U$ is free. Their quadratic phase is determined by this factor, while the critical decay threshold also depends on the positive factor.
    \item Finally, the propagators of \eqref{intro.Schrodinger} are metaplectic operators. Therefore, Theorem \ref{thm-main3} implies dynamical versions of Hardy's uncertainty principle when applied to concrete Schr\"odinger evolutions. Specifically, using Theorem \ref{thm-main3} we obtain results in the following form: under the notation of \eqref{intro.Schrodinger}, 
    \begin{equation}
       \text{If $u_0\in L^p(\bR)$ and $\exists t_0>0$ s.t. $|u(t_0,x)|\lesssim e^{-\pi \gamma x^2}$ with $\gamma> \gamma_0=\gamma_0(t_0)$}\quad \Longrightarrow\quad u=0,
    \end{equation}
    where the critical exponent $\gamma_0$ is computed explicitly and it might be included or not in the inequality above, according to the structure of the corresponding Hamiltonian flow.
    We refer to Section \ref{sec:Examples} for the exact statements for the heat equation, the Davies equation and the complex harmonic oscillator.
    \end{enumerate} 

    Altogether, our results show that the metaplectic semigroup setting gives rise to a richer and more fragile extremal theory, where the existence and form of the extremizers depend on subtle algebraic conditions on the entries of the symplectic matrix, going beyond what the classical Hardy uncertainty principle would suggest. 
    
    \medskip
    
    {\bf Outline.} This work is structured as follows. The notation and the prerequisites are summarized in Section \ref{sec:prerequisites}. Importantly, we detail the standard Phragmén-Lindelöf argument for proving Hardy's uncertainty principle in Section \ref{standardArgument}.
    In Section \ref{sec:positivity}, we prove Theorem \ref{intro.Main-1}. In Section \ref{sec:IRP}, we prove the integral representations in Theorem \ref{intro.main-2}. Hardy's uncertainty principle is proven in Section \ref{secHUP}, together with Gaussian bounds for $|\widehat Sf|$ and continuity of $\widehat Sf$. Finally, dynamic versions of Hardy's uncertainty principle are formulated in Section \ref{sec:Examples}. We conclude in Section \ref{sec:conclusions} with a discussion of higher-dimensional extensions and perturbed propagators.
    
\section{Prerequisites}\label{sec:prerequisites}
    \subsection{Notation and matrix algebra}
    Let us summarize the main notation used in this work.
    
    \noindent We denote by $\bF$ the field of real or complex numbers and by $\bF^d$ the corresponding finite dimensional spaces. We will consider the standard inner products in these spaces,  $x\cdot y=x^\top y$, $x,y\in\rd$,  and $x\cdot y=x^\top\overline{y}$, $x,y\in\cd$.
    We write $\bF^{m\times n}$ for the space of $m\times n$ matrices with entries in $\bF$. If $S\in\bF^{d\times d}$, $S^\top$ and $S^\ast$ denote the transpose and the conjugate-transpose of $S$, respectively, and we shall write $S\geq0$ if $S$ is positive semi-definite, i.e., $S=S^*$ and $z^\ast Sz\geq0$ for every $z\in\bF^{d}$. As it shall not cause confusion, the same notation is used to denote the adjoint of operators.
    
    If $G$ is a semigroup with identity  element, we denote by $G^\times$ the group of units of $G$, that is the group of the invertible elements of $G$.
     
     We denote by $\mathbb{U}$ the (closed) upper-half plane of $\mathbb{C}$, whereas $\bT=\{z\in\bC:|z|=1\}$ is the unit circle. Its elements are called {\em phase factors}.
     We denote by $\Sym(d,\bF)$ the class of symmetric matrices with entries in $\bF$. 
    
    \subsection{Symplectic group}
     This subsection is devoted to symplectic analysis. We refer the reader to \cite{CorderoBook,deGosson,GRT2026,Hormander3} as main references.
    Let 
    \begin{equation}\label{defJ}
        J=\begin{pmatrix}
            0 & 1\\
            -1 & 0
        \end{pmatrix}
    \end{equation} 
    be the matrix of the canonical symplectic form of $\bR^2$. A matrix $S\in\bC^{2\times 2}$ is {\em symplectic} if $S^\top JS=J$. By writing
    \begin{equation}\label{blockS}
        S=\begin{pmatrix}
            a & b\\
            c & d
        \end{pmatrix}, \qquad a,b,c,d\in\bC,
    \end{equation}
    this is equivalent to the determinant condition
    \begin{equation}\label{blockSrel}
        \det(S)=ad-bc=1.
    \end{equation}
    
    We denote by $\Sp(\bC)$ the complex symplectic group. 
    We shall write the real and the imaginary parts of the entries of $S$ with subscripts: $a=a_1+ia_2$, and similarly for $b$, $c$ and $d$. We extend this notation to $S$, and denote $S_1=\Re(S)$ and $S_2=\Im(S)$.
    If the entries of $U\in\Sp(\bC)$ are real, we write $U\in\Sp(\bR)$. Matrices $S\in\Sp(\bC)$ with blocks \eqref{blockS} having $b\neq0$ are called {\em free}.

    The symplectic group $\Sp(\bF)$, with $\bF=\bR$ or $\bF=\bC$ is generated by $J$, and by matrices in the subgroups $\{\cD_\lambda\}_{\lambda\neq0}$ and $\{V_\alpha\}_{\alpha\in\bF}$, where
    \begin{equation}\label{defDE}
        \cD_\lambda=\begin{pmatrix}
            \lambda^{-1} & 0\\
            0 & \lambda
        \end{pmatrix},\qquad
   {\text{and}}\qquad 
        V_\alpha =\begin{pmatrix}
            1 & 0\\
            \alpha & 1
        \end{pmatrix}.
    \end{equation}

    The group $\Sp(\bC)$ has a subsemigroup $\Sp_+(\bC)$ of matrices satisfying the {\em positivity} condition
    \begin{equation}\label{posCond}
        i^{-1}\big(S^\ast JS-J\big)\geq0.
    \end{equation}
    \begin{definition}
        The semigroup $\Sp_+(\bC)$ is the semigroup of {\em positive} $2\times2$ symplectic matrices.
    \end{definition}
    The semigroup $\Sp_+(\bC)$ was investigated by Brunet and Kramer \cite{BK1,BK2}, Olshanskii \cite{Olshanskii}, and H\"ormander \cite{Hormander1}.
    \begin{remark}\label{Rem23}
        By Sylvester's criterion for positive semi-definiteness of Hermitian matrices, see e.g., \cite[Chapter 7.6]{Meyer}, we can rephrase \eqref{posCond} as follows. First, observe that
        \begin{equation}
            i^{-1}\big(S^\ast JS-J\big)=\begin{pmatrix}
                2\Im(a^\ast c) & i^{-1}(a^\ast d-c^\ast b-1)\\
                i(d^\ast a-b^\ast c-1) & 2\Im(b^\ast d)
            \end{pmatrix}.
        \end{equation}
        Observe that it is also $i^{-1}\big(S^\ast JS-J\big)=-2\Im(S)^\top JS$, which implies \begin{equation}
            \det(\Im(S))={\frac{1}{4}}\det\big[i^{-1}\big(S^\ast JS-J\big)\big]\geq0.
        \end{equation}
        
        By Sylvester's criterion, we see that if $S\in\Sp(\bC)$, then \eqref{posCond} holds if and only if the principal minors are non-negative, i.e.,
        \begin{align}
            \label{Sstar1}
           & c_2a_1-a_2c_1\geq0\\
           \label{Sstar2}
           & d_2b_1-b_2d_1\geq0,\\
           \label{Sstar3}
           & \det(\Im(S))\geq0.
        \end{align}
       Together with the condition $ad-bc=1$, \eqref{Sstar3} implies that $\det(\Re(S))\geq1$.
       
    \end{remark}

    The following result is proven in \cite{GRT2026}. 
    \begin{proposition}\label{propGenSppiu}
        A set of generators for the semigroup $\Sp_+(\bC)$ is 
        \begin{equation}
            \{J\}\cup \{\cD_\lambda\}_{\lambda\in\bR^\times}\cup\{V_\alpha\}_{\alpha\in\mathbb{U}}.
        \end{equation}
        \end{proposition}
 \begin{remark}\label{remDE}
        It is worth observing that matrices $\cD_\lambda$ as in \eqref{defDE} belong to $\Sp_+(\bC)$ if and only if $\lambda\in\bR^\times$, see e.g. \cite{GRT2026}.
    \end{remark}
    \subsection{The metaplectic semigroup}
    The content of this section synthesizes the results in \cite{GiacchiAchar2026} and the classical theory in \cite{Brunet,Hormander3,Howe,tramontana}.
    Consider the complexification of the Schr\"odinger representation of the Heisenberg group $\bH$:
    \begin{equation}\label{SchroRep}
        \rho(x,\xi;\tau)f(y)=e^{i\tau}e^{-i\pi x\xi}e^{2\pi i\xi y}\int_{-\infty}^\infty \hat f(\eta)e^{2\pi i(y-x)\eta}\mathrm d\eta, \qquad y\in\bR,\; x,\xi\in\bC,\; \tau\in\bC.
    \end{equation}
    initially defined on a suitable dense subspace of $L^2(\bR)$. If $x,\xi,\tau\in\bR$, the operator \eqref{SchroRep} is unitary on $L^2(\bR)$ and $\rho:\bH\to\mathcal U(L^2(\bR))$ is a unitary irreducible representation of $\bH$.
    If $S\in\Sp_+(\bC)$ there exists a non-zero continuous operator $\widehat S:L^2(\bR)\to L^2(\bR)$ such that
    \begin{equation}\label{intert}
        \widehat S\rho(z;\tau)=\rho(Sz;\tau)\widehat S, \qquad z\in\bR^2,\;\tau\in\bR
    \end{equation}
    holds for every $f\in L^2(\bR)$ as an identity between $L^2$ functions.
    This intertwiner $\widehat S$ is never unique, but if $\widehat S'$ is another continuous non-zero operator satisfying \eqref{intert}, then $\widehat S'=c\widehat S$ for a constant $c\in\bC\setminus\{0\}$.
    The family of {\em metaplectic operators}, i.e.,
    \begin{equation}
        \{\widehat S:S\in{\Sp}_+(\bC)\},
    \end{equation}
    forms a semigroup under composition, which admits a subsemigroup of contractions consisting precisely of two operators for each $S\in\Sp_+(\bC)$. This subsemigroup is the {\em oscillator semigroup} (or {\em metaplectic semigroup}), denoted here by $\Mp_+(\bC)$. It is the two-fold cover of $\Sp_+(\bC)$, and the {\em projection} $\pi_+^{Mp}:\widehat S\in\Mp_+(\bC)\mapsto S\in\Sp_+(\bC)$ is a semigroup homomorphism with kernel $\{\pm\mathrm{id}_{L^2}\}$.
    For fixed $S\in\Sp_+(\bC)$, the properties of any lift $\widehat S$ relevant to our uncertainty principles are invariant under multiplication by a non-zero scalar. We shall therefore omit non-zero scalar normalization factors throughout the presentation of our results.
    
    The group of units of $\Mp_+(\bC)$ is the {\em metaplectic group}, which admits the following integral representation: every $U\in\Sp(\bR)$ with entries as in \eqref{blockS} lifts to a pair of operators with integral representation
    \begin{equation}\label{defMuS}
        \widehat U f(x)=\begin{cases}
            |b|^{-1/2}e^{i\pi dx^2/b}\int_{\bR} f(y)e^{i\pi ay^2/b}e^{-2\pi ixy/b}\ud y, & \text{if $U$ is free},\\
            |a|^{-1/2}e^{i\pi cx^2/a}f(x/a), & \text{otherwise,}
        \end{cases}\qquad f\in\cS(\bR).
    \end{equation}
    We call {\em quadratic Fourier transforms} the operators $\widehat U$ corresponding to free symplectic matrices. Observe that such operators have kernels
    \begin{equation}\label{defKU}
        K_U(x,y)=|b|^{-1/2}e^{2\pi iq(x,y)}, \qquad q(x,y)=\frac 1 2 \cQ(x,y)\cdot(x,y), \, x,y\in\bR,
    \end{equation}
    where
    \begin{equation}
        \cQ=\begin{pmatrix}
            d/b & -1/b\\
            -1/b & a/b
        \end{pmatrix}.
    \end{equation}
    \begin{example}\label{FractionalFT}
    Important instances of these operators are {\em fractional Fourier transforms}, defined for $\theta\in\bR\setminus\{k\pi:k\in\mathbb Z\}$ by
    \begin{equation}\label{defFracFT}
        \Fur_\theta f(x)=i^{-1/2}|\sin\theta|^{-1/2}e^{i\pi \cot\theta x^2}\int_{-\infty}^\infty f(y)e^{i\pi\cot\theta y^2}e^{-2\pi ixy/\sin\theta}\mathrm dy,\qquad f\in\cS(\bR).
    \end{equation}
    and for $\theta=k\pi$ ($k\in\bZ$) by $\Fur_\theta f(y)=f((-1)^ky)$.
    The projection of $\Fur_\theta$ is
    \begin{equation}
        \Pi_\theta=\begin{pmatrix}
            \cos\theta & \sin\theta\\
            -\sin\theta & \cos\theta
        \end{pmatrix}.
    \end{equation}
    \end{example}
    To distinguish metaplectic operators in $\Mp_+(\bC)$ from those in $\Mp(\bR)$, we refer to the latter as {\em unitary metaplectic operators}, the reason is the following proposition.
    \begin{proposition}
    The following continuity properties hold.
          \begin{enumerate}[(i)]
              \item Every operator $\widehat S\in\Mp_+(\bC)$ is injective, and it is surjective (and unitary) if and only if $\widehat S\in\Mp(\bR)$.
              \item The inverses of the operators in $\Mp_+(\bC)$ are closed, densely defined on $L^2(\bR)$. Moreover, $\Mp_+(\bC)$ is closed under adjunction.
              \item The $L^2$-contractions $\widehat S\in \Mp_+(\bC)$ restrict to continuous operators from $\cS(\bR)$ to itself, and they extend to continuous injective operators from $\cS'(\bR)$ to itself by duality:
              \begin{equation}
                  \la \widehat Sf,g\ra =\la f,\widehat S^\ast g\ra, \qquad f\in\cS'(\bR),\; g\in\cS(\bR),
              \end{equation}
              where $\widehat S^\ast$ is the adjoint of $\widehat S$.
          \end{enumerate}
      \end{proposition}
  By lifting the generators of $\Sp_+(\bC)$ listed in Proposition \ref{propGenSppiu}, a neat description of the generators of $\Mp_+(\bC)$ was obtained in \cite{GRT2026}.
  \begin{proposition}\label{propGenerators}
      The semigroup $\Mp_+(\bC)$ is generated by operators in the form
      \begin{align}
        \label{defFT}
          & \Fur f(\xi)=i^{-1/2}\widehat f(\xi)=i^{-1/2}\int_{\bR}f(y)e^{-2\pi i\xi y}\ud y, \qquad f\in\cS(\bR),\\
          \label{defTE}
          &  {\mathfrak{T}_\lambda f(x)=i^{m}|\lambda|^{1/2}f(\lambda x)}, 
          \qquad \lambda\in\bR^\times,\\
          \label{defpa}
          &\mathfrak p_{\alpha}f(x)=e^{i\pi\alpha x^2}f(x), \qquad \alpha\in\mathbb{U},
      \end{align}
      where $m\in\bZ$ satisfies $(-1)^m=\mathrm{sign}(\lambda)$ ({\em Maslov index}).
  \end{proposition}
  As aforementioned, the phase factors $i^{-1/2}$ and $i^m$, and the implicit sign ambiguities, appearing in \eqref{defFT} and \eqref{defTE} do not play an active role in the present analysis, and they will be therefore omitted for simplicity.

  We conclude this section with an integral representation formula for the metaplectic operators associated with $U\in\Sp(\bR)$ having entry $a\in\bR^\times$. These operators play a fundamental role in expressing the propagators of quadratic Schr\"odinger evolutions.
    \begin{proposition}
        Let $\widehat U\in\Mp(\bR)$ have projection $U$ with entries \eqref{blockS}. If $a\in\bR^\times$, then
        \begin{equation}\label{seeA}
            \widehat Uf(x)=|a|^{-1/2}e^{i\pi cx^2/a}\int_{\bR}\widehat f(\eta)e^{-i\pi b\eta^2/a}e^{2\pi i\eta x/a}\ud\eta, \qquad f\in\cS(\bR).
        \end{equation}
    \end{proposition}

    \subsection{The Gaussian integral and Phragmén-Lindelöf principle}\label{standardArgument}
    We consider $L^2(\bR)$ endowed with the sesquilinear inner product $\la f,g\ra=\int_{\bR} f(y)\overline{g(y)}\ud y$. This inner product extends uniquely to a sesquilinear duality pairing $\cS'(\bR)\times\cS(\bR)$, where $\cS(\bR)$ denotes the Schwartz class of rapidly decaying complex-valued smooth functions on $\bR$, and $\cS'(\bR)$ is the space of tempered distributions. 
    
    We shall also use the following result concerning the Fourier transform of the Gaussian integrals.
    \begin{theorem}\cite[Appendix A]{Folland}\label{thm:Gaussian}
        Let $\lambda\in\bC^\times$ have $\Re(\lambda)> 0$. Then, for every $\zeta\in\bC$,
        \begin{equation}
            \int_{\bR}e^{-\pi\lambda y^2}e^{-2\pi i\zeta y}\ud y=\lambda^{-1/2} e^{-\pi\lambda^{-1}\zeta^2},
        \end{equation}
        where $\lambda^{-1/2}=|\lambda|^{-1/2}e^{-\frac i2 \arg(\lambda)}$, 
        $-\frac\pi 2\le \arg \lambda\le \frac \pi 2.$
    \end{theorem}
    \begin{remark}
If $\lambda\in i\bR\setminus\{0\}$, the function
$e^{-\pi\lambda y^2}$ defines a tempered distribution, and
\[
    \Fur\big(e^{-\pi\lambda y^2}\big)(\zeta)
    =\lambda^{-1/2}e^{-\pi\lambda^{-1}\zeta^2}
    \qquad\text{in\  }\ \cS'(\bR),
\]
where $\lambda^{-1/2}$ is defined by continuous extension
of the chosen branch from the right half-plane.
This identity follows by replacing $\lambda$ with
$\lambda+\varepsilon$ in the Gaussian integral formula above 
and letting $\varepsilon\downarrow0$.
\end{remark}
    In this work, we consider operators on $L^2(\bR)$. If  $T:L^2(\bR)\to L^2(\bR)$ is linear and bounded, we denote by $T^\ast$ its adjoint, and $T$ is self-adjoint if $T=T^\ast$, whereas it is positive if $\la Tf,f\ra\geq 0$ for every $f\in L^2(\bR)\setminus\{0\}$. Recall that in this context, every positive operator is also self-adjoint.

We recall the Phragm\'en--Lindel\"of principle in a sector
\cite[Chapter~VI, Corollary~4.2]{Conway}.

\begin{lemma}\label{lem:PL}
Let $\Gamma\subset\bC$ be an open sector with vertex at the
origin and opening $\alpha\in(0,2\pi)$. Suppose that $g$ is
holomorphic in $\Gamma$ and that, for some $M\geq0$,
\[
    \limsup_{\Gamma\ni z\to\zeta}|g(z)|\leq M
    \qquad\text{for every }\zeta\in\partial\Gamma.
\]
If there exist $A,B>0$ and $0<\sigma<\pi/\alpha$ such that
\[
    |g(z)|\leq A e^{B|z|^\sigma},
    \qquad z\in\Gamma,
\]
then $|g(z)|\leq M$ in $\Gamma$.
\end{lemma}

\begin{corollary}\label{cor:PL-critical}
Let $\phi$ be an entire function. Suppose that, for some
$A,C>0$,
\[
    |\phi(x+iy)|\leq C e^{Ax^2},
    \qquad
    |\phi(x)|\leq C,
    \qquad x,y\in\bR.
\]
Then $\phi$ is constant.
\end{corollary}

\begin{proof}
Fix $0<\theta<\pi/2$ and consider
\[
    \phi(z)\exp\left(\frac{iA}{2}\cot\theta\,z^2\right)
\]
in the sector $0<\arg z<\theta$. On the positive real axis
the exponential factor has modulus one. On the ray
$\arg z=\theta$, we have $y=x\tan\theta$, so its modulus is
$e^{-Axy\cot\theta}=e^{-Ax^2}$. Thus the function is bounded
by $C$ on both boundary rays. In the sector its modulus is at
most $Ce^{A|z|^2}$. Since $2<\pi/\theta$, Lemma~\ref{lem:PL}
gives the bound $C$ throughout the sector.

For each $z$ in the first quadrant, letting $\theta\uparrow\pi/2$
gives $|\phi(z)|\leq C$. Applying the same argument to
$\phi(-z)$, $\overline{\phi(\bar z)}$ and
$\overline{\phi(-\bar z)}$ gives the bound in the remaining
quadrants. Liouville's theorem now implies that $\phi$ is constant.
\end{proof}

  \subsection{An operator perspective on the positivity condition}\label{sec:positivity}
  
  We derive integral representations for operators in $\Mp_+(\bC)$ directly in terms of the entries of their positive symplectic projections
  We begin by proving Theorem \ref{intro.Main-1}.

  \begin{proof}[Proof of Theorem \ref{intro.Main-1}]
Write
\[
    S=\begin{pmatrix}a&b\\c&d\end{pmatrix},
    \qquad ad-bc=1,
\]
and set
\[
    H_S=i^{-1}(S^\ast JS-J).
\]
The matrix $H_S$ is Hermitian, and by \eqref{posCond},
$S\in\Sp_+(\bC)$ if and only if $H_S\geq0$.

\noindent
\emph{(i) The case $b=0$.}
This case is included in \cite{GRT2026}; we give a direct
verification for completeness.
Since $ad=1$, we have $a\neq0$ and $d=a^{-1}$.
A direct computation gives
\[
    (H_S)_{11}=2\Im(\overline{a}c),
    \qquad
    (H_S)_{12}=i^{-1}(\overline{a}a^{-1}-1),
    \qquad
    (H_S)_{22}=0.
\]
A Hermitian $2\times2$ matrix with vanishing lower-right
entry is positive semidefinite if and only if its
off-diagonal entries vanish and its upper-left entry
is nonnegative. Therefore
\[
    H_S\geq0
    \quad\Longleftrightarrow\quad
    a=\overline{a}
    \ \text{ and }\ 
    \Im(\overline{a}c)\geq0.
\]
The first condition is equivalent to $a\in\bR^\times$.
For such $a$,
\[
    \Im(\overline{a}c)=a^2\Im(c/a),
\]
so the second condition is equivalent to
$\Im(c/a)\geq0$. This proves $(i)$.

\noindent
\emph{(ii) The case $b\neq0$.}
Define
\[
    Q=
    \begin{pmatrix}
        d/b & -1/b\\
        -1/b & a/b
    \end{pmatrix},
    \qquad
    L=
    \begin{pmatrix}
        a & b\\
        1 & 0
    \end{pmatrix}.
\]
We claim that
\begin{equation}\label{eq:positivity-congruence}
    H_S=2L^\ast(\Im Q)L.
\end{equation}
Here $\Im Q$ denotes the entrywise imaginary part of $Q$.
Since $Q$ is symmetric, we have
\[
    Q-Q^\ast=2i\,\Im Q.
\]
Moreover, using $ad-bc=1$, we obtain
\[
    QL=
    \begin{pmatrix}
        c & d\\
        0 & -1
    \end{pmatrix},
    \qquad
    L^\ast QL=
    \begin{pmatrix}
        \overline{a}c & \overline{a}d-1\\
        \overline{b}c & \overline{b}d
    \end{pmatrix}.
\]
Subtracting the adjoint of the latter matrix gives
\[
    L^\ast QL-(L^\ast QL)^\ast=S^\ast JS-J.
\]
Therefore,
\[
    2iL^\ast(\Im Q)L
    =L^\ast(Q-Q^\ast)L
    =S^\ast JS-J,
\]
which proves \eqref{eq:positivity-congruence}.

Finally, $\det L=-b\neq0$. Since congruence by an
invertible matrix preserves positive semidefiniteness,
$S\in\Sp_+(\bC)$ if and only if $
    \Im Q\geq0.$
This is precisely \eqref{intro.posDefQF}, and proves $(ii)$.
\end{proof}
  
    A first corollary is the following.
    \begin{corollary}\label{corbReal}
        Let $S\in\Sp_+(\bC)$ have entries \eqref{blockS} with $b\neq 0$. If $\Imab=0$, then $b\in\bR^\times$. In particular, if $\Im(b)\neq0$, then $\Imab>0$.
    \end{corollary}
    \begin{proof}
        It follows by \eqref{intro.posDefQF} which, if $\Im(a/b)=0$, implies
      \begin{equation}
          0\leq \det\begin{pmatrix}
              \Im(d/b) & -\Im(1/b)\\
              -\Im(1/b) & 0
          \end{pmatrix}=-\Im(1/b)^2,
      \end{equation}
      that is $\Im(1/b)=0$, and we have done.
      
    \end{proof}
    
  As a corollary, we obtain a characterization of matrices in $\Sp_+(\bC)$ with invertible top-left entry.
  \begin{corollary}\label{corollaryA}
      Let $S\in\Sp(\bC)$ have entries as in \eqref{blockS},
with $a\neq0$. Then, $S\in\Sp_+(\bC)$ if and only if
      \begin{equation}\label{posDefQF-A}
          \Im\begin{pmatrix}
              c/a & 1/a\\
              1/a & -b/a
          \end{pmatrix}\geq0.
      \end{equation}
  \end{corollary}
  \begin{proof}
      It follows by Theorem \ref{intro.Main-1}, applied to
      \begin{equation}\label{defSprime}
         S'= SJ^\top=
          \begin{pmatrix}
              b & -a\\
              d & -c
          \end{pmatrix},
      \end{equation}
      whose top-right entry $-a$ is invertible.
  \end{proof}

    \begin{remark}
        The positivity conditions \eqref{intro.posDefQF} and \eqref{posDefQF-A} are not themselves surprising. In the real setting, the quadratic forms defined through the matrices appearing therein are related to the integral kernels of the corresponding metaplectic operators, see  \eqref{defKU} and \eqref{seeA}.
    \end{remark}

\section{Integral representation formulas}\label{sec:IRP}
  In this section, we use the construction of the metaplectic semigroup in \cite{Hormander3} to obtain explicit integral representation formulas for metaplectic operators, in terms of the corresponding symplectic projections. Ultimately, we will prove Theorem \ref{intro.main-2}, but let us stress a few facts before.
  \begin{remark}
      The operator in \eqref{into.freeMetapOps} has kernel
      \begin{equation}
          K_S(x,y)=|b|^{-1/2}e^{2\pi iq(x,y)}, \qquad q(x,y)=\frac 1 2 \cQ(x,y)\cdot(x,y), \qquad x,y\in\bR,
      \end{equation}
      where
      \begin{equation}
          \cQ=\begin{pmatrix}
              d/b & -1/b\\
              -1/b & a/b
          \end{pmatrix}\in\Sym(2,\bC)
      \end{equation}
      has $\Im(\cQ)\geq0$ by Theorem \ref{intro.Main-1}. In particular, $K_S\in\cS'(\bR^2)$ and $\widehat S$ is well-defined on $\cS(\bR)$.
  \end{remark}
    \begin{remark}
      Observe that if $b\neq0$, then
      \begin{equation}\label{decompSb}
          S=V_{d/b}\cD_{1/b}JV_{a/b},
      \end{equation}
      so one may be tempted to infer
      \begin{equation}\label{decompSb2}
          \widehat S=\mathfrak{p}_{d/b}\mathfrak{T}_{1/b}\Fur \mathfrak p_{a/b},
      \end{equation}
      and then use this factorization to obtain a representation formula by composition.
      However, as observed in Remark \ref{remDE}, the matrix $\cD_{1/b}$ is in $\Sp_+(\bC)$ if and only if $b\in\bR^\times$, and consequently, $\mathfrak{T}_{1/b}\in\Mp_+(\bC)$ if and only if $b\in \bR^\times$. 
      Therefore, the operator $\mathfrak{T}_{1/b}$ is not defined if $b\notin\bR$.
      That is to say, the factorization in \eqref{decompSb}, which holds in $\Sp(\bC)$, does not hold in $\Sp_+(\bC)$ for general $b\in\bC^\times$, and this reflects on the fact that the corresponding factorization \eqref{decompSb2} does not make sense unless $b$ is real. 
  \end{remark}
  \begin{remark}
      If $\Im(a/b)>0$, formula \eqref{into.freeMetapOps} holds in $L^\infty(\bR)$. If $f(x)$ has a Gaussian decay, the integral in \eqref{into.freeMetapOps} converges also when $\Im(a/b)=0$.
  \end{remark}
 
    \subsection{Proof of Theorem \ref{intro.main-2}}
        The case $b=0$ is proven in \cite{GRT2026}, so we focus on the case where $b\neq0$. In this case, the quadratic form
        \begin{equation}
            \cQ=\begin{pmatrix}
                d/b & -1/b\\
                -1/b & a/b
            \end{pmatrix}
        \end{equation}
        has positive semi-definite imaginary part and, consequently, 
        \begin{equation}\label{kernel}
            K_S(x,y)=|b|^{-1/2}e^{2\pi iq(x,y)}, \qquad q(x,y)=\frac 1 2 \cQ(x,y)\cdot(x,y), \, x,y\in\bR,
        \end{equation}
        defines a tempered distribution in $\cS'(\bR^2)$. We must prove that the operator
        \begin{equation}
            Tf(x)=\int_{\bR}f(y)K_S(x,y)\ud y, \qquad f\in\cS(\bR)
        \end{equation}
        coincides, up to a phase factor, with $\widehat S$.

        The canonical relation determined by $q$ is
$
    \Lambda_q
    =\{(x,\partial_xq,y,-\partial_yq):x,y\in\bC\}.
$
We denote
\[
    \eta=-\partial_yq=\frac{x-ay}{b},
\]
and obtain $x=ay+b\eta$ and, since $ad-bc=1$,
$
    \partial_xq=\frac{dx-y}{b}=cy+d\eta.
$
Therefore
\[
    \Lambda_q
    =\{(ay+b\eta,cy+d\eta,y,\eta):y,\eta\in\bC\}
    =\{(Sz,z):z\in\bC^2\},
\]
which is the graph of $S$.

We refer to the Gaussian-kernel construction in
\cite[Section 5]{Hormander3}. Then the kernel of
$\widehat S$ is
\[
    \gamma_S e^{2\pi iq(x,y)},
    \qquad
    \gamma_S^2=i\,\partial_x\partial_yq=-\frac{i}{b},
\]
where the choice of the square root corresponds to
the choice of the lift $\widehat S$.
Since $|\gamma_S|=|b|^{-1/2}$, this kernel coincides
with $K_S$ up to a phase. Thus
$\widehat Sf=\zeta Tf$ for some $\zeta\in\bT$,
which concludes the proof.

   \subsection{An alternative representation formula} 
   Similarly to Corollary \ref{corollaryA}, we obtain an integral representation formula for the case $a\neq0$. When $a,b\neq0$ these formulas represent the same operator differently.
  \begin{corollary}\label{CorollaryA2}
      Let $\widehat S\in\Mp_+(\bC)$ have projection $S$ with entries \eqref{blockS}. If $a\neq0$, then
      \begin{equation}\label{inrRepSa}
          \widehat Sf(x)=|a|^{-1/2}e^{i\pi cx^2/a}\int_{\bR}\widehat f(\eta)e^{-i\pi b\eta^2/a}e^{2\pi i\eta x/a}\ud\eta, \qquad f\in\cS(\bR),
      \end{equation}
      up to a phase.
  \end{corollary}
  \begin{proof}
      If $a\neq0$, then 
      \begin{equation}
          S=S'J,
      \end{equation}
      where $S'$ is the matrix in \eqref{defSprime}. The assertion follows then by choosing corresponding metaplectic lifts
      \begin{equation}\label{defSprimeA2}
          \widehat Sf=\widehat S'\widehat f,
      \end{equation}
      up to a phase, and writing the integral expression of $\widehat S'$ given by Theorem \ref{intro.main-2}.
      
  \end{proof}
    
\section{Proof of the main result}\label{secHUP}
In this section, we show Hardy's uncertainty principle for operators in $\Mp_+(\bC)$. Before delving into the details, let us observe that if the projection $S\in\Sp_+(\bC)$ of $\widehat S\in\Mp_+(\bC)$ has entries \eqref{blockS} with $b=0$, then $\widehat Sf$ preserves the compactness of supports, exactly as in the case of unitary metaplectic operators. Indeed, in this case, $a\in\bR^\times$ and it follows by Corollary \ref{CorollaryA2} that
\begin{equation}
    \widehat Sf(x)=|a|^{-1/2}e^{i\pi cx^2/a}f\big(x/a\big),
\end{equation}
up to a phase. For this reason, in what follows, we shall assume $b\neq0$. 

\subsection{Intrinsic Gaussian decay and continuity}
Mapping properties of integral operators with complex Gaussian kernels are classical, see \cite[Propositions 2.1, 3.1 and 4.1]{GN2004}, and \cite{Lieb} for sharp norm results.
We record the following elementary estimate to make the intrinsic Gaussian decay rate explicit in terms of the entries of the symplectic matrix.

\begin{proposition}\label{propFreeDecay}
    Let $\widehat S\in\Mp_+(\bC)$ have projection $S$ with entries \eqref{blockS} and $b\neq0$. If $\Im(a/b)>0$, then
    \begin{equation}\label{positivityOfBound}
        \Im(d/b)-\Im(a/b)^{-1}\Im(b)^2/|b|^4\geq0
    \end{equation}
    and
    \begin{equation}\label{freeDecay}
        |\widehat Sf(x)|\leq |b|^{-1/2}\big(p'\Imab\big)^{-\frac{1}{2p'}}\norm{f}_p e^{-\pi\big(\Im(d/b)-\Im(a/b)^{-1}\Im(b)^2/|b|^4\big)x^2},
    \end{equation}
    for every $f\in L^p(\bR)$, $1\leq p\leq\infty$, where we understand $\big(p'\Imab\big)^{-\frac{1}{2p'}}=1$ for $p=1$. Moreover, the integral in \eqref{into.freeMetapOps} defines a continuous representative of $\widehat Sf$.
\end{proposition}
\begin{proof}
    The positivity \eqref{positivityOfBound} follows by Theorem \ref{intro.Main-1} and the identity
    \begin{equation}
        \Im(d/b)-\Im(a/b)^{-1}b_2^2/|b|^4=\Im(a/b)^{-1}\det\left(\Im\begin{pmatrix}
            d/b & -1/b\\
            -1/b & a/b
        \end{pmatrix}\right).
    \end{equation}
    By Theorem \ref{intro.main-2} $(ii)$,
    \begin{align}
        |\widehat Sf(x)|\leq|b|^{-1/2} e^{-\pi\Im(d/b)x^2}\int_{\bR}|f(y)|e^{-\pi\Im(a/b)y^2}e^{-2\pi b_2xy/|b|^2}\ud y.
    \end{align}
    If $1<p\leq\infty$, by H\"older's inequality and Theorem \ref{thm:Gaussian},
    \begin{align}
        |\widehat Sf(x)|&\leq |b|^{-1/2}\norm{f}_pe^{-\pi\Im(d/b)x^2}\Big(\int_{\bR}e^{-\pi p'\Im(a/b)y^2}e^{-2\pi b_2p'xy/|b|^2}\ud y\Big)^{1/p'}\\
        &=|b|^{-1/2}\big(p'\Imab\big)^{-\frac{1}{2p'}}\norm{f}_pe^{-\pi\big(\Im(d/b)-\Im(a/b)^{-1}b_2^2/|b|^4\big)x^2}.
    \end{align}
    If $p=1$, then
    \begin{align}
        |\widehat Sf(x)|&\leq |b|^{-1/2}\norm{f}_1e^{-\pi\Im(d/b)x^2}{\sup}_{y\in\bR}e^{-\pi\Im(a/b)y^2}e^{-2\pi b_2xy/|b|^2}.
    \end{align}
    For every $x\in\bR$, the function
    \begin{equation}
        h_x(y)=e^{-\pi\Im(a/b)y^2}e^{-2\pi b_2xy/|b|^2}
    \end{equation}
    attains its maximum at $y_M=-\Imab^{-1}b_2x/|b|^2$. Since
    \begin{equation}
        h_x(y_M)=e^{\pi\Im(a/b)^{-1}b_2^2x^2/|b|^4}
    \end{equation}
    the estimate follows also for $p=1$. Finally, for every $R>0$ and $|x|\leq R$,
    \[
        |f(y)e^{i\pi ay^2/b}e^{-2\pi ixy/b}| \leq C_R |f(y)|e^{-\pi\Imab y^2/2}.
    \]
    The right-hand side is integrable by H\"older's inequality, including the endpoints $p=1$ and $p=\infty$. Dominated convergence in \eqref{into.freeMetapOps} therefore proves the continuity of $\widehat Sf$, thereby concluding the proof.

\end{proof}

We stress that Proposition \ref{propFreeDecay} entails a natural decay condition in the dissipative framework: if $\Im(a/b)>0$ and $f\in L^\infty(\bR)$, then 
\begin{equation}\label{decaySfx}
    |\widehat Sf(x)|\lesssim e^{-\pi\big(\Im(d/b)-\Im(a/b)^{-1}\Im(b)^2/|b|^4\big)x^2}, \qquad x\in\bR.
\end{equation}
This explains the apparently technical decay condition \eqref{intro.decaySf1}.

\subsection{A prototypical example}\label{ExampleHeat}
Before proving Theorem \ref{thm-main3}, we present an example that heuristically illustrates the possible expected cases arising in the proof of Theorem \ref{thm-main3}.

Let us consider a metaplectic operator associated with
    \begin{equation}
        S=\begin{pmatrix}
            1 & b\\
            0 & 1
        \end{pmatrix}, \qquad \Im(b)<0.
    \end{equation}
    Up to a phase,
    \begin{equation}\label{exampleb}
        \widehat Sf(x)=\Fur^{-1}\big(e^{-i\pi b(\cdot)^2}\widehat f\big)(x).
    \end{equation}
    Operators of this type form a subsemigroup of $\Mp_+(\bC)$ that includes the heat semigroup. We assume for simplicity $|b|=1$.
    Assume that
    \begin{equation}
        |\widehat Sf(x)|\lesssim e^{-\pi\gamma x^2}, \qquad x\in\bR
    \end{equation}
    (observe that this is exactly \eqref{intro.decaySf1} in this particular case). We proceed heuristically by applying the Fourier transform to $\widehat Sf$:
    \begin{equation}\label{ineqEx1}
        |\Fur\widehat Sf(\xi)|\lesssim \norm{\widehat f}_\infty e^{\pi b_2\xi^2}.
    \end{equation}
    Assume that $ \widehat f\in L^\infty(\bR)$, so that the right-hand side of \eqref{ineqEx1} is finite for every $\xi\in\bR$. 
    By Theorem \ref{intro.HardyClassical}, if $-\gamma b_2>1$, then $f=0$, whereas for $\gamma=-1/b_2$ (recall that $b_2<0$), we have
    \begin{equation}\label{exampleSf}
       \widehat S f(x)=C' e^{\pi x^2/b_2}.
    \end{equation}
    To determine $f$, we use \eqref{exampleb} and write
    \begin{equation}
        \Fur^{-1}\big(e^{-i\pi b(\cdot)^2}\widehat f\big)(x)=C' e^{\pi x^2/b_2},
    \end{equation}
    we obtain
    \begin{equation}
        \widehat f(\xi)=C e^{i\pi b_1\xi^2}.
    \end{equation}
    If $b_1=0$, then $f(x)=C\delta$, so if we further assume $f\in L^\infty(\bR)$, we expect 
    \[
    -b_2\gamma=1 \qquad \Longrightarrow \qquad f=0.
    \]
    Instead, if $b_1\neq0$,
    \[
    -b_2\gamma=1 \qquad \Longrightarrow \qquad f(x)=C e^{-i\pi x^2/b_1},
    \]
    which is now admissible under the assumption $f\in L^\infty(\bR)$. Observe that if we require $f\in L^1(\bR)$ instead, then $f=0$ whenever $-b_2\gamma=1$.
    
We are ready to prove Theorem \ref{thm-main3}.
\subsection{Proof of Theorem \ref{thm-main3}}
    In a line, Item $(i)$ follows from classical Hardy's uncertainty principle, whereas Item $(ii)$ requires deeper analytic and geometrical considerations.
    \begin{itemize}
        \item[$(i)$] If $\Imab=0$, then $b\in\bR^\times$, and we can write
    \begin{align}
     e^{-i\pi dx^2/b}\widehat Sf(x)
    =\zeta |b|^{-1/2}\widehat g(x/b),
       \qquad x\in\bR,
    \end{align}
    where $g(x)=f(x)e^{i\pi ax^2/b}$ and $\zeta\in\bT$. By \eqref{intro.decaySf2} we have
    \begin{align}
        |g(x)|=|f(x)|\lesssim e^{-\pi \delta x^2},\qquad x\in\bR,
    \end{align}
    whereas
    \begin{align}
        |\widehat g(\xi)|
    =|b|^{1/2}
      \left|e^{-i\pi db\xi^2}\widehat Sf(b\xi)\right|
    \lesssim e^{-\pi b^2\gamma\xi^2},
    \qquad \xi\in\bR.
    \end{align}
    Therefore, the assertion follows by applying Theorem \ref{intro.HardyClassical} to $g$.
    \item[$(ii)$] If $\Im(a/b)>0$, we argue using the Phragmén-Lindel\"of principle, see Lemma \ref{lem:PL}.
    The rest of this proof is divided into two steps.\\
    {\em Step 1. Case $|b|=1$.}
    In this case, the decay assumption \eqref{intro.decaySf1} simplifies as
    \begin{equation}\label{decaySf1b1}
        |\widehat Sf(x)|\lesssim e^{-\pi\big(\gamma+\Im(d/b)-\Im(a/b)^{-1}b_2^2\big)x^2}, \qquad x\in\bR.
    \end{equation}
    Define
    \begin{equation}\label{defg}
        g(z)=\widehat Sf(z)=e^{i\pi dz^2/b}\int_{\bR}f(y)e^{i\pi ay^2/b}e^{-2\pi izy/b}\ud y, \qquad z\in\bC
    \end{equation}
    and $G(z)=e^{-i\pi dz^2/b}g(z)$. 
    By differentiating under the integral, we obtain
    \begin{equation}\label{Holomorph}
         \partial_zG(z)=-\frac{2\pi i}{b}\int_{\bR}yf(y)e^{i\pi ay^2/b}e^{-2\pi izy/b}\ud y.
     \end{equation}
     Since $\Im(a/b)>0$, we have
\[
    \left|y e^{i\pi ay^2/b}e^{-2\pi izy/b}\right|
    =|y|e^{-\pi\Im(a/b)y^2+2\pi y\Im(z/b)}.
\]
For $z$ in a compact subset of $\bC$, this is bounded by
\[
    C|y|e^{-\frac{\pi}{2}\Im(a/b)y^2}.
\]
After multiplication by $|f(y)|$, this bound is
integrable by H\"older's inequality. The same argument
applies to the integrand defining $G$, without the
factor $y$. Thus differentiation under the integral
is justified, locally uniformly in $z$, and $G$ is entire.
     
    Theorem \ref{thm:Gaussian} yields for $p\neq1$,
         \begin{align}\label{estimateGx}
             |G(z_1+iz_2)|&\lesssim\norm{f}_p\Big(\int_{\bR}e^{-\pi p'\Im(a/b)y^2}e^{2\pi p'y(z_2b_1-z_1b_2)}\ud y
             \Big)^{1/p'}\\
             &=\norm{f}_p\big(p'\Imab\big)^{-\frac{1}{2p'}} e^{\pi\Im(a/b)^{-1}(z_2b_1-z_1b_2)^2}.
         \end{align}
         For $p=1$ the upper bound \eqref{estimateGx} for $|G(z)|$ follows (with constant $1$) by observing that
         \begin{equation}
             e^{-\pi \Im(a/b)y^2}e^{2\pi y(z_2b_1-z_1b_2)}\leq e^{\pi\Im(a/b)^{-1}(z_2b_1-z_1b_2)^2}.
         \end{equation}
         Observe that \eqref{estimateGx} can be rephrased as
         \begin{equation}
             |G(z)|\lesssim  e^{\pi\Im(a/b)^{-1}\Im(z/b)^2},\qquad z\in\bC,
         \end{equation}
         which implies the quadratic exponential bound $|G(z)|\lesssim e^{\pi\Im(a/b)^{-1}|z|^2}$.
        Then, the function 
    \begin{equation}\label{defphi}
        \phi(z):=G(z)e^{\pi \Im(a/b)^{-1}b_1 z^2/b}, \qquad z\in\bC
    \end{equation}
    satisfies for $z=x+iy\in\bC$,
    \begin{align}
    |\phi(z)|
    &\lesssim
    e^{\frac{\pi}{\Im(a/b)}
    \left[
        (b_1y-b_2x)^2
        +\Re\left(\frac{b_1z^2}{b}\right)
    \right]}
    =e^{\frac{\pi}{\Im(a/b)}
        (b_1^2+b_2^2)x^2}
    =e^{\pi x^2/\Im(a/b)},
\end{align}
we have used that $|b|=1$.
    By \eqref{decaySf1b1}, we can also estimate
    \begin{align}
        |\phi(x)|&=|\widehat Sf(x)|e^{\pi\Im(d/b)x^2}e^{\pi \Im(a/b)^{-1}b_1^2x^2}\\
        &\lesssim e^{-\pi\big(\gamma+\Im(d/b)-\Im(a/b)^{-1}b_2^2\big)x^2}e^{\pi\Im(d/b)x^2} e^{\pi \Im(a/b)^{-1}b_1^2x^2}\\
        &=e^{-\pi\big(\gamma-\Im(a/b)^{-1}\big)x^2}.
    \end{align}
    Consequently, if $\gamma\Im(a/b)\geq1$, there exists $C>0$
such that
\[
    |\phi(x+iy)|\leq C e^{\pi x^2/\Im(a/b)},
    \qquad
    |\phi(x)|\leq C,
    \qquad x,y\in\bR.
\]
The preceding two
estimates and Corollary~\ref{cor:PL-critical} imply that
$\phi$ is constant.
        \begin{enumerate}[(a)]
        \item {\em The super-critical regime.} If $\gamma\Im(a/b)>1$, the real-axis estimate gives
$\phi(x)\to0$ as $x\to+\infty$, so the constant function $\phi$
equals zero. 
        This entails $\widehat Sf=0$, whence $f=0$ by injectivity. 
        \item {\em The critical regime.} If $\gamma\Im(a/b)=1$,  $\phi(x)={C'}$ for some $C'\in\bC$.  Hence 
        \begin{equation}\label{limit}
            \widehat Sf(x)={C'} e^{i\pi dx^2/b}e^{-\pi\Im(a/b)^{-1}b_1x^2/b}.
        \end{equation}
        If $C'=0$, then $f=0$ by injectivity and there is nothing to prove. We may therefore assume that $C'\neq0$. 
        Observe that by \eqref{into.freeMetapOps},
        \begin{equation}
            \widehat S\delta_0(x)=C_0e^{i\pi dx^2/b}, \qquad C_0\neq0,
        \end{equation}
        which coincides with \eqref{limit} when $b_1=0$. Therefore, if $b_1=0$, the injectivity of $\widehat S$ on $\cS'$ forces $f=(C'/C_0)\delta_0$. This is an $L^p$ function if and only if $C'=0$. 
        Assume now that $b_1\neq0$ and set
        \begin{equation}\label{deff0}
            f_0(x)=e^{-i\pi(a_1/b_1)x^2}.
        \end{equation}
        Since
        \begin{equation}
            \frac ab-\frac{a_1}{b_1}=i\frac{\Imab}{bb_1},
        \end{equation}
        Gaussian integration and \eqref{into.freeMetapOps} yield
        \begin{align}
            \widehat Sf_0(x)=C_1e^{i\pi dx^2/b}e^{-\pi \Imab^{-1}b_1x^2/b}
        \end{align}
        for some $C_1\neq0$. By injectivity, it must be $f=(C'/C_1)f_0$. This function is in $L^p$ if and only if $p=\infty$. 
    \item {\em The sub-critical regime.} 
    Assume $0\leq\gamma\Im(a/b)<1$. We first show that there exists $\lambda\in\mathbb U$ with $\Im\lambda>0$ such that $f(x)=e^{i\pi\lambda x^2}$ satisfies \eqref{decaySf1b1}. 
    For $\Im(\lambda)>0$, Gaussian integration, applied to \eqref{into.freeMetapOps}, and the identity
    \[
        \frac{c+d\lambda}{a+b\lambda}=\frac db-\frac{1}{b(a+b\lambda)},
    \]
    which follows from $ad-bc=1$, give
    \begin{align}
        |\widehat Sf(x)|&=|a+b\lambda|^{-1/2}e^{-\pi\Im(d/b)x^2}e^{\pi \Im\big(((a+b\lambda)^{-1})/b\big)x^2}\\
        &=|a+b\lambda|^{-1/2}e^{-\pi\big(\Im(d/b)-\Im(a/b)^{-1}b_2^2\big)x^2}\\
        &\qquad\times\; e^{-\pi\big(\Im(a/b)^{-1}b_2^2-\Im(((a+b\lambda)^{-1})/b)\big)x^2}.
    \end{align}
    We wonder for which $\gamma\Im(a/b)<1$ there is $\lambda\in\mathbb{U}$ such that
    \begin{align}
        e^{-\pi\big(\Im(d/b)-\Im(a/b)^{-1}\Im(b)^2\big)x^2}&e^{-\pi\big(\Im(a/b)^{-1}\Im(b)^2-\Im(((a+b\lambda)^{-1})/b)\big)x^2}\\
        &\lesssim e^{-\pi\big(\gamma+\Im(d/b)-\Im(a/b)^{-1}\Im(b)^2\big)x^2},
    \end{align}
    or equivalently
    \begin{equation}\label{systemLimit}
            \gamma\leq \Im(a/b)^{-1}{b_2^2}-\Im\Big(\frac{1}{b(a+\lambda b)}\Big).
    \end{equation}
    Put $\eta=\Im(a/b)>0$ and recall that $|b|=1$.
If $b_1\neq0$, take
\[
    \lambda_\varepsilon=-\frac{a_1}{b_1}+i\varepsilon,
    \qquad \varepsilon>0.
\]
Then
\[
    \lim_{\varepsilon\downarrow0}
    \left[
        \frac{b_2^2}{\eta}
        -\Im\left(\frac{1}{b(a+b\lambda_\varepsilon)}\right)
    \right]
    =\frac{b_2^2+b_1^2}{\eta}
    =\frac1\eta.
\]
If $b_1=0$, take
\[
    \lambda_R=-\Re(a/b)+iR,\qquad R>0.
\]
Since $b^2=-1$, we have
\[
    \frac{b_2^2}{\eta}
    -\Im\left(\frac{1}{b(a+b\lambda_R)}\right)
    =\frac1\eta-\frac{1}{\eta+R}
    \longrightarrow\frac1\eta
    \qquad\text{as }R\to\infty.
\]
Thus, for every $\gamma<\eta^{-1}$, there exists
$\lambda$ with $\Im\lambda>0$ such that
\eqref{systemLimit} holds strictly. The corresponding
Gaussian belongs to every $L^p(\bR)$,
$1\leq p\leq\infty$, and satisfies \eqref{decaySf1b1}.
    
    It remains to check that in the sub-critical regime, the subspace
    \begin{equation}
        \Omega=\{f\in L^2(\bR) : \eqref{decaySf1b1} \; \text{holds}\}
    \end{equation}
    is dense in $L^2(\bR)$.
    Integration by parts shows that
    \begin{equation}
        \widehat S\big((\cdot)e^{i\pi\lambda (\cdot)^2}\big)(x)=\frac{x}{a+b\lambda}\widehat Sf(x),
    \end{equation}
    and a standard induction argument proves that for every $k$,
    \begin{equation}
        \widehat S\big((\cdot)^ke^{i\pi\lambda (\cdot)^2}\big)(x)=p(x)\widehat S\big(e^{i\pi\lambda (\cdot)^2}\big)(x),
    \end{equation}
    where $p$ is a polynomial of degree at most $k$ with complex coefficients. So, if $P$ is a polynomial of degree $k$, then
    \begin{align}
        \big|\widehat S \big(P e^{i\pi\lambda (\cdot)^2}\big)(x)\big|\lesssim _{\varepsilon,P} e^{\pi\varepsilon x^2}\big|\widehat S\big(e^{i\pi \lambda (\cdot)^2}\big)(x)\big|
    \end{align}
    for every $\varepsilon>0$. 
    Let $\gamma<\Imab^{-1}$ and $\varepsilon>0$ sufficiently small so that $\gamma'=\gamma+\varepsilon<\Imab^{-1}$. Then, there is $\lambda\in\mathbb U$, $\lambda_2>0$, so that $|\widehat S(e^{i\pi\lambda(\cdot)^2})(x)|\lesssim e^{-\pi(\gamma'+\Im(d/b)-\Imab^{-1}b_2^2)x^2}$.
    For such a $\lambda$, we have
    \begin{align}
        |\widehat S(P e^{i\pi\lambda(\cdot)^2})(x)|&\lesssim_{\varepsilon,P} e^{\pi\varepsilon x^2}e^{-\pi\big(\gamma'+\Im(d/b)-\Imab^{-1}b_2^2\big) x^2}\\
        &=e^{-\pi\big(\gamma+\Im(d/b)-\Imab^{-1}b_2^2\big) x^2}.
    \end{align}
    We have proven that for every $\gamma<\Imab^{-1}$ the functions $P(x)e^{i\pi \lambda x^2}$ satisfy \eqref{decaySf1b1} for a suitable choice of $\lambda$.
    For this fixed $\lambda$ with $\Im\lambda>0$, the span
of the functions $P(x)e^{i\pi\lambda x^2}$ is dense
in $L^2(\bR)$. Indeed, this span is obtained from
the span of the Hermite functions by a unitary
dilation and multiplication by
$e^{i\pi\Re(\lambda)x^2}$.

    \noindent
    {\em Step 2. Case $|b|\neq1$.}
    The positive symplectic matrix
    \begin{equation}
        S':=S\begin{pmatrix}
            |b| & 0\\
            0 & |b|^{-1}
        \end{pmatrix},
    \end{equation}
    has entries
    \begin{equation}
        S'=\begin{pmatrix}
            a|b| & b/|b|\\
            c|b| & d/|b|
        \end{pmatrix}=:\begin{pmatrix}
            a' & b'\\
            c' & d'
        \end{pmatrix}, \qquad |b'|=1.
    \end{equation}
    Observe that the assumption $b'\neq0$ is still satisfied by $S'$, whereas $\Im(a'/b')=|b|^2\Imab$ and $\Im(d'/b')=\Im(d/b)$. 
    From the operator perspective, 
    \begin{equation}
        \widehat S=\widehat S'\mathfrak T_{|b|},
    \end{equation}
    up to a phase.
    By defining $g=|b|^{1/2}f(|b|\cdot)\in L^p(\bR)$, we have $\widehat Sf=\widehat S'g$. Moreover, the decay assumption \eqref{intro.decaySf1} implies 
    \begin{align}
        |\widehat S'g(x)|&=|\widehat Sf(x)|\lesssim e^{-\pi\big(\gamma+\Im(d/b)-\Im(a/b)^{-1}b_2^2/|b|^4\big)x^2}\\
        &=e^{-\pi\big(\gamma+\Im(d'/b')-\Im(a'/b')^{-1}b_2'^2\big)x^2}.
    \end{align}
    So $|\widehat S'g|$ satisfies the corresponding decay assumption \eqref{decaySf1b1} of {\em Step 1}. 
    The theorem then follows for $\widehat{S}$ after observing that
    \[
        \gamma \Im(a'/b') = \gamma \Imab |b|^2,
    \]
    and that the change of variables $\mathfrak{T}_{|b|}$ is an isomorphism on $L^p(\mathbb{R})$ for every $1\leq p\leq\infty$.
    
 \end{enumerate} 
 \end{itemize}
 The proof is now concluded.

\subsection{Hardy's uncertainty principle for \texorpdfstring{$\Fur L^p$}{}}
A similar uncertainty principle for $\Fur L^p$ spaces follows now as a corollary of Theorem \ref{thm-main3} and Corollary \ref{CorollaryA2}.

\begin{theorem}\label{thm-main4}
    Let $\widehat S\in\Mp_+(\bC)$ have projection $S$ as in \eqref{blockS} with $a,b\neq0$. Assume $\Im(b/a)<0$.
        Suppose that $f\in \Fur L^p(\bR)$ for some $1\leq p\leq\infty$ and
        \begin{equation}\label{decaySfxA}
            |\widehat Sf(x)|\lesssim e^{-\pi\big(\gamma+\Im(c/a)+\Im(b/a)^{-1}\Im(a)^2/|a|^4\big)x^2}, \qquad x\in\bR,
        \end{equation}
        for some $\gamma\geq 0$.
        \begin{itemize}
            \item If $-\gamma \Im(b/a)> |a|^{-2}$, then $f=0$.
            \item If $-\gamma \Im(b/a)=|a|^{-2}$ 
            then $f=0$, unless $p=\infty$
            
            and $\Re(a)\neq0$.
In the latter case, for some $C\in\bC$:
\begin{itemize}
    \item[-] If $\Re(b)=0$, then $f=C\delta_0$.
    \item[-] If $\Re(b)\neq0$, then
    \begin{equation}\label{limitCaseA2}
        f(x)=C e^{-i\pi\Re(a)x^2/\Re(b)}.
    \end{equation}
\end{itemize}
     \item If $-\gamma \Im(b/a)<|a|^{-2}$, there is a dense subspace of $L^2(\bR)$ satisfying \eqref{decaySfxA}.
        \end{itemize}
\end{theorem}
\begin{proof}
    Write $\widehat S=\widehat S'\Fur$, as in the proof of Corollary \ref{CorollaryA2}, see \eqref{defSprimeA2}. 
    Then, the result follows by applying Theorem \ref{thm-main3} to $\widehat S'$ and $\widehat f\in L^p(\bR)$. 
    In particular, for $-\gamma\Im(b/a)=|a|^{-2}$, we have $f=0$ unless $p=\infty$ and $a_1\neq0$. In the latter case, 
    \begin{equation}
        \widehat f(\xi)=C'e^{i\pi b_1\xi^2/a_1}.
    \end{equation}
    If $b_1=0$, then $f=C\delta_0$ for some $C\in\bC$. Otherwise, we obtain \eqref{limitCaseA2} by Fourier inversion.
    
\end{proof}

\begin{proof}[Proof of Theorem \ref{intro.subcriticalHermite}]
Define 
$F(z)=e^{\pi\omega_S z^2}\widehat Sf(z).$
The integral representation \eqref{into.freeMetapOps} and the assumption $\Im(a/b)>0$
show that $F$ is entire. By H\"older's inequality
and Gaussian integration, for $z=x+iy$,
\begin{align}
    |F(z)|
    \lesssim
    \exp\left\{
        \frac{\pi}{\Im(a/b)}\bigl(\Im(z/b)\bigr)^2
        +\pi\Re\left(
            -\frac{i\gamma_0\Im b}{b}z^2
        \right)
    \right\}
    =e^{\pi\gamma_0 y^2}.
\end{align}
We used that
$\omega_S+id/b=-i\gamma_0\Im(b)/b$ and
$\gamma_0=(\Im(a/b)|b|^2)^{-1}$.
Moreover,
\[
    \Re\omega_S
    =\Im(d/b)-\frac{(\Im b)^2}{\Im(a/b)|b|^4},
\]
so the decay assumption \eqref{intro.decaySf1} gives
\[
    |F(x)|\lesssim e^{-\pi\gamma x^2},
    \qquad x\in\bR.
\]

Set
\[
    \rho=\sqrt{\gamma/\gamma_0}\in(0,1),
    \qquad
    \mu=\frac{1-\rho}{1+\rho}.
\]
Since $h(z)=\ell^{1/2}F(\ell z)$ and
$\ell=(\gamma\gamma_0)^{-1/4}$, we obtain
\[
    |h(x)|\lesssim e^{-\pi\rho x^2},
    \qquad
    |h(x+iy)|\lesssim e^{\pi y^2/\rho}.
\]
In particular, $h\in L^2(\bR)$.

Following the argument of Vemuri in \cite{Vemuri2008}, we consider the Bargmann transform of $h$,
\[
    \mathcal Bh(w)
    =2^{1/4}\int_{\bR}
        h(t)e^{-\pi t^2+2\sqrt{\pi}\,tw-w^2/2}
        \ud t.
\]
With our normalization of the Hermite functions,
\(
    \mathcal Bh(w)
    =\sum_{n=0}^{\infty}
        \langle h,\psi_n\rangle
        \frac{w^n}{\sqrt{n!}}.
\)

The real-axis bound for $h$ and Gaussian integration
give, for $w=u+iv$,
\[
    |\mathcal Bh(u+iv)|
    \lesssim e^{(\mu u^2+v^2)/2}.
\]

To obtain a second bound on $|\mathcal Bh(u+iv)|$, we shift the contour of
integration from $\R$ to $\bR+is$, where
\[
    s=\frac{\rho v}{\sqrt{\pi}(1+\rho)}.
\]
The bound for $h(x+iy)$ implies that the integrals
over the vertical segments joining $\pm x$ to
$\pm x+is$ tend to zero as $x\to\infty$. More precisely,
for fixed $w=u+iv$ and $s$,
\[
    \left|
    \int_{\pm x}^{\pm x+is}
        h(t)e^{-\pi t^2+2\sqrt{\pi}\,tw-w^2/2}
        \ud t
    \right|
    \lesssim_{w,s}
    e^{-\pi x^2+2\sqrt{\pi}|u|x}
    \longrightarrow 0.
\]
Cauchy's theorem therefore allows us to shift the
contour to $\bR+is$.
Using the same bound on the shifted contour yields
\begin{align}
    |\mathcal Bh(u+iv)|
    \lesssim
    \exp\left\{
        \frac{u^2+v^2}{2}
        +\pi(1+\rho^{-1})s^2
        -2\sqrt{\pi}\,sv
    \right\}
    =e^{(u^2+\mu v^2)/2}.
\end{align}

These are precisely the two bounds used in
\cite[proofs of Theorems 2.1 and 2.2]{Vemuri2008},
after writing its Bargmann variable as $\sqrt{2}\,w$. We provide the details for the sake of the reader.
We apply the Phragm\'en--Lindel\"of principle to
$e^{i\sqrt{\mu}\,w^2/2}\mathcal Bh(w)$
in the sector
\[
    \arctan\sqrt{\mu}
    \leq \arg w
    \leq \frac{\pi}{2}-\arctan\sqrt{\mu}.
\]
Indeed, for $w=re^{i\theta}$, the two preceding
estimates give
\begin{align*}
    \left|
        e^{i\sqrt{\mu}\,w^2/2}\mathcal Bh(w)
    \right|
    \lesssim
    \min\left\{
        e^{\frac{r^2}{2}
            (\sqrt{\mu}\cos\theta-\sin\theta)^2},
        e^{\frac{r^2}{2}
            (\cos\theta-\sqrt{\mu}\sin\theta)^2}
    \right\}.
\end{align*}
The first estimate gives boundedness on the ray
$\arg w=\arctan\sqrt{\mu}$, and the second gives
boundedness on the ray
$\arg w=\pi/2-\arctan\sqrt{\mu}$.

Since $0<\mu<1$, the opening of the sector is
$\pi/2-2\arctan\sqrt{\mu}<\pi/2$. The function
$e^{i\sqrt{\mu}\,w^2/2}\mathcal Bh(w)$ has growth bounded
by $Ce^{C|w|^2}$ and is bounded on both boundary rays.
Hence Lemma~\ref{lem:PL}, applied with $\sigma=2$, gives
boundedness throughout the sector.
Consequently,
\begin{equation}
    |\mathcal Bh(re^{i\theta})|
    \lesssim
    \exp\left\{
        \frac{\sqrt{\mu}}{2}r^2\sin(2\theta)
    \right\},
    \qquad
    \arctan\sqrt{\mu}
    \leq\theta\leq
    \frac{\pi}{2}-\arctan\sqrt{\mu}.
    \label{eq:Hermite-sector-bound}
\end{equation}

For $0\leq\theta\leq\arctan\sqrt{\mu}$ and $\frac{\pi}{2}-\arctan\sqrt{\mu}\le \theta\le\pi/2$, we have \[\min\{\mu\cos^2\theta+\sin^2\theta,\cos^2\theta+\mu\sin^2\theta\}\leq 2\mu/(1+\mu).\]
Therefore the original estimates imply
\[
    |\mathcal Bh(re^{i\theta})|
    \lesssim \exp\left(\frac{\mu}{1+\mu}r^2\right).
\]
We can repeat the argument in each of the quadrants.

We now integrate over the circle of radius $r\geq1$.
In the sector of \eqref{eq:Hermite-sector-bound}, write
$\theta=\pi/4+t$. Since
$\sin(2\theta)=\cos(2t)\leq 1-8t^2/\pi^2$ for
$|t|\leq\pi/4$, we obtain
\begin{align*}
    \int_0^{2\pi} |\mathcal Bh(re^{i\theta})|\,\ud\theta
    \lesssim e^{\mu r^2/(1+\mu)}
      +e^{\sqrt{\mu}\,r^2/2}
       \int_{-\pi/4}^{\pi/4}
       e^{-4\sqrt{\mu}\,r^2t^2/\pi^2}\,\ud t \lesssim_{\mu} r^{-1}e^{\sqrt{\mu}\,r^2/2},
\end{align*}
where we used $\mu/(1+\mu)<\sqrt{\mu}/2$.

By Cauchy's formula and the Hermite expansion of $\mathcal Bh$,
\[
    |\langle h,\psi_n\rangle|
    \lesssim_{\mu}
    \sqrt{n!}\,r^{-n-1}e^{\sqrt{\mu}\,r^2/2}.
\]
For $n\geq1$, take $r^2=n/\sqrt{\mu}$. Stirling's formula
then implies
\(
    |\langle h,\psi_n\rangle|
    \lesssim_{\mu} n^{-1/4}\mu^{n/4}.
\) This concludes the proof.
\end{proof}
   \begin{remark}
Under the assumptions of Theorem~\ref{intro.subcriticalHermite},
the normalized output $h$ satisfies a Gaussian bound under every
fractional Fourier transform, see \eqref{defFracFT} above for the definition. More precisely,
\[
    \sup_{\theta\in\bR} |\Fur_\theta h(x)|
    \lesssim
    \exp\left(
        -\pi\frac{\sqrt{\gamma}}
        {\sqrt{\gamma_0}+\sqrt{\gamma_0-\gamma}}\,x^2
    \right),
    \qquad x\in\bR.
\]
Indeed, choose $\alpha>0$ so that
$\tanh(2\alpha)=\sqrt{\gamma/\gamma_0}$.
Then Theorem~\ref{intro.subcriticalHermite} gives
\[
    |\langle h,\psi_n\rangle|
    \lesssim (n+1)^{-1/4}e^{-\alpha n}.
\]
Since $\Fur_\theta$ acts on each $\psi_n$ by a scalar of modulus
one, the weighted Hermite estimate of Radchenko and Ramos
\cite{RadchenkoRamos2025} yields
\[
    |\Fur_\theta h(x)|
    \lesssim \sum_{n=0}^{\infty}
        (n+1)^{-1/4}e^{-\alpha n}|\psi_n(x)|
    \lesssim e^{-\pi\tanh(\alpha)x^2},
\]
uniformly in $\theta$. The stated bound follows from
$\tanh(\alpha)=\sqrt{\gamma}/
(\sqrt{\gamma_0}+\sqrt{\gamma_0-\gamma})$.
\end{remark}

    \section{Dynamic versions of Hardy's uncertainty principle}\label{sec:Examples}
    
    In this section, we apply the theory developed so far to evolution problems in the form
    \begin{equation}\label{CauchyEvo}
        \begin{cases}
           \frac{i}{2\pi}\partial_tu=a^\w(x,D) u, \qquad \text{on $(0,+\infty)\times\bR$},\\
            u(0,\cdot)=u_0\in L^p(\bR),
        \end{cases}
    \end{equation}
    where $a^\w(x,D)$ is the Weyl quantization of a complex
quadratic form $a$ whose imaginary part is negative semidefinite.
We write
\begin{equation}
    a(x,\xi)=\frac12 e\xi^2+fx\xi-\frac12 gx^2,
    \qquad e,f,g\in\bC,
\end{equation}
with $\Im a(x,\xi)\leq0$ for all $(x,\xi)\in\bR^2$.
For $u_0\in\cS(\bR)$, the Weyl quantization is given by
  \begin{equation}
a^\w(x,D)u_0(x)=\int_{\bR^2}a\Big(\frac{x+y}{2},\xi\Big)u_0(y)e^{2\pi i(x-y)\xi}\ud y\ud \xi,
    \end{equation}
    and the integral is interpreted in the distributional sense.
    
    Hereafter, we shall consider $1\leq p\leq \infty$ and apply Theorem \ref{thm-main3} to infer dynamical uncertainty inequalities for specific PDEs. Similar conclusions hold if we consider $u_0\in \Fur L^p(\bR)$, by invoking Theorem \ref{thm-main4} instead, with the necessary precautions on the critical scenarios.

    \subsection{The heat equation}
    Let us show that the uncertainty principle for the heat equation in Theorem \ref{intro.ThmMalinnikova} and Section \ref{ExampleHeat} is now a straightforward consequence of our results.
    \begin{theorem}\label{thm:heat}
        Let $e^{-2\pi ita^\w(x,D) }$ be the propagator of \eqref{heat}, with $u_0\in L^p(\bR)$, $1\leq p\le\infty$. If there exists $t_0>0$ such that $|u(t_0,x)|\lesssim e^{-\pi\gamma x^2}$ for some $\gamma\geq(4\pi t_0)^{-1}$, then $u=0$.
    \end{theorem}
    \begin{proof}
        The Hamiltonian flow of the heat equation \eqref{heat} is
        \begin{equation}
            S_t=\begin{pmatrix}
                1 & -4\pi it\\
                0 & 1
            \end{pmatrix}.
        \end{equation}
        Observe that $\Im(-1/(4\pi it))=1/(4\pi t)>0$ and $\Re(-4\pi it)=0$ for every $t>0$.
        By applying Theorem \ref{thm-main3} we obtain the assertion.
        
    \end{proof}

    \subsection{The non-self-adjoint harmonic oscillator}
    Let us now consider the non-self-adjoint harmonic oscillator, or Davies equation:
    \begin{equation}\label{Davies}
            \frac i {2\pi}\partial_tu=a^\w(x,D) u, \qquad u(0,x)=u_0\in L^p(\bR),\qquad 1\leq p\leq\infty,
    \end{equation}
    where 
    \begin{equation}
        a^\w(x,D) =i\left(\frac{e^{-i\vartheta}}{8\pi^2}\partial_x^2-\frac 1 2e^{i\vartheta}x^2\right),\qquad |\vartheta|<\pi/2.
    \end{equation}
    We establish Hardy's uncertainty principle for the evolution of \eqref{Davies}.
    
    \begin{theorem}\label{thm:Davies}
        Consider problem \eqref{Davies}. If there exists $t_0>0$ such that
        \begin{equation}
            |u(t_0,x)|\lesssim e^{-\pi\big(\gamma+\tanh(t_0)\cos(\vartheta)\big)x^2}, \qquad x\in\bR
        \end{equation}
        for some $\gamma\geq\frac{2}{\sinh(2t_0)\cos(\vartheta)}$, then $u=0$, with only one exception: if $\gamma=\frac{2}{\sinh(2t_0)\cos(\vartheta)}$, $\vartheta\neq0$ and $p=\infty$, then
        \begin{equation}\label{DaviesExtremizer}
            u_0(x)=C e^{i\pi\coth(t_0)x^2/\sin(\vartheta)}.
        \end{equation}
    \end{theorem}
    \begin{proof}
    The Hamiltonian flow can be computed explicitly, see \cite{Viola}:
    \begin{equation}
        S_t=\begin{pmatrix}
            \cosh(t) & -i\sinh(t)e^{-i\vartheta}\\
            i\sinh(t) e^{i\vartheta} & \cosh(t)
        \end{pmatrix}=:\begin{pmatrix}
            a(t) & b(t)\\
            c(t) & d(t)
        \end{pmatrix}, \qquad t>0.
    \end{equation}
    Observe that the condition $|\vartheta|<\pi/2$ guarantees that $S_t\in\Sp_+(\bC)$ for every $t>0$. This can be easily verified by means of the condition \eqref{posDefQF-A}, that in this case reads as
    \begin{equation}
        \begin{pmatrix}
            \tanh(t)\cos(\vartheta) & 0\\
            0 &\tanh(t)\cos(\vartheta)
        \end{pmatrix}>0.
    \end{equation}
    Moreover, we remark that $\Im\big(a(t)/b(t)\big)=\cos(\vartheta)/\tanh(t)>0$ and 
    $\Re\big(b(t)\big)=-\sinh(t)\sin(\vartheta)\neq0$ for every $t>0$ and $-\pi/2<\vartheta<\pi/2$, $\vartheta\neq0$.
    It follows then by Theorem \ref{thm-main3} that if $\tanh(t_0)\cos(\vartheta)\gamma>\cosh(t_0)^{-2}$ for $t=t_0$, then, $u_0=0$.  
    In terms of $\gamma$, the critical condition
$\gamma\Im\big(a(t_0)/b(t_0)\big)=|b(t_0)|^{-2}$ reads as
    \begin{align}
        \gamma\,\frac{\cos(\vartheta)}{\tanh(t_0)}=\frac{1}{\sinh(t_0)^{2}}
        \qquad\Longleftrightarrow\qquad
        \gamma=\frac{1}{\sinh(t_0)\cosh(t_0)\cos(\vartheta)}
        =\frac{2}{\sinh(2t_0)\cos(\vartheta)}.
    \end{align}
    Moreover, if $\gamma=2/(\sinh(2t_0)\cos(\vartheta))$, then also $u=0$ unless $\vartheta\neq0$ and $p=\infty$. 
    It remains to identify the extremizers in this exceptional case. By Theorem
    \ref{thm-main3} $(ii)$, they are given by \eqref{intro.def-fb20}, so we
    evaluate the relevant quantities at $t=t_0$. From the expressions of $a(t_0)$ and $b(t_0)$, we infer
    \begin{align}
        &\Re\big(a(t_0)\big)=\cosh(t_0),\\
        &\Re\big(b(t_0)\big)=-\sinh(t_0)\sin(\vartheta).
    \end{align}
    Consequently,
    \begin{align}
        -\frac{\Re\big(a(t_0)\big)}{\Re\big(b(t_0)\big)}
        =\frac{\coth(t_0)}{\sin(\vartheta)},
    \end{align}
    and \eqref{intro.def-fb20} reads as \eqref{DaviesExtremizer}.
    This concludes the proof. 
    
    \end{proof}

    \begin{remark}
Theorem \ref{thm:Davies} provides a concrete instance of the breakdown
of the classical extremal picture discussed in the Introduction.
For fixed $t_0>0$ and $\vartheta\neq0$, the critical extremizers
in Theorem~\ref{thm:Davies} are precisely the scalar multiples
of the quadratic chirp in \eqref{DaviesExtremizer}. As in the
classical Hardy theorem, they form a one-dimensional space.
Here, however, every nonzero extremizer has constant modulus,
so it belongs to $L^\infty(\bR)$ and to no $L^p(\bR)$ with
$1\leq p<\infty$.
This is in sharp contrast with Theorems
\ref{intro.HardyClassical} and \ref{intro.thm13}.

This exceptional family degenerates as $\vartheta\to0$. Indeed, the
coefficient of the quadratic phase in \eqref{DaviesExtremizer} becomes
unbounded, so that the extremizers oscillate on increasingly fine
scales. At the limiting value $\vartheta=0$, the Hamiltonian flow
satisfies $\Re b(t)=0$, and the exceptional critical family ceases to
exist: the zero function is then the only admissible datum at the
critical threshold.
\end{remark}
    
\subsection{Complex harmonic oscillator}\label{subsec:CHO}
    We now consider the complex harmonic oscillator
    \begin{equation}\label{CHO}
        \frac i {2\pi}\partial_tu=\zeta\left(\frac {1}{8\pi^2}\partial_x^2
        -\frac {1} 2x^2\right)u, \qquad u(0,\cdot)=u_0\in L^p(\bR), \qquad 1\leq p\leq \infty,
    \end{equation}
    where $\zeta=\zeta_1+i\zeta_2, \zeta_1,\zeta_2\in\bR$ and $\zeta_2>0$. 
   Hardy's uncertainty principle is now a matter of applying either Theorem \ref{thm-main3} or Theorem \ref{thm-main4}.

   \begin{theorem}
       Under the notation of this subsection, 
       \begin{align}\label{propEx3}
            u(t,x)=|\sin(\zeta t)|^{-1/2}e^{-i\pi \cot(\zeta t) x^2}\int_{\bR}u_0(y)e^{-i\pi \cot(\zeta t)y^2}e^{2\pi ixy/\sin(\zeta t)}\ud y,
        \end{align}
        $t>0$ and $x\in\bR$, up to a phase.
       Moreover, if there exists $t_0>0$ such that
       \begin{equation}\label{decaySfCHO}
           |u(t_0,x)|\lesssim e^{-\pi\big(\gamma+\tanh(\zeta_2t_0)\big)x^2}
       \end{equation}
       where $\gamma\geq 2/\sinh(2\zeta_2t_0)$,
       then $u=0$, with the following exception: if $u_0\in L^\infty(\bR)$, $\sin(\zeta_1t_0)\neq0$ and $\gamma=2/\sinh(2\zeta_2t_0)$, in which case 
            \begin{equation}
                u_0(x)=Ce^{i\pi \cot(\zeta_1t_0)x^2
}.
        \end{equation}
   \end{theorem}
    \begin{proof}
        For the sake of completeness, we derive the integral expression of the propagator of \eqref{CHO}. This is usually obtained by means of generalized Mehler formulas. Here, we use directly Theorem \ref{intro.main-2}.
        
        The Hamiltonian flow associated with \eqref{CHO} is $S_t=e^{tX}$, where
        \begin{align}
            X=J\begin{pmatrix}
                -\zeta & 0\\
                0 & -\zeta
            \end{pmatrix}=-\zeta J,
        \end{align}
        so that
        \begin{equation}
            X^{2k}=(-1)^k\zeta^{2k}I \qquad \text{and} \qquad X^{2k+1}=(-1)^{k+1}\zeta^{2k+1}J, \qquad k=0,1,2,\ldots.
        \end{equation}
        Then,
        \begin{align}\label{CHOSt}
            S_t&=\begin{pmatrix}
                \cos(\zeta t) & -\sin(\zeta t)\\
                \sin(\zeta t) & \cos(\zeta t)
            \end{pmatrix}=:\begin{pmatrix}
                a(t) & b(t)\\
                c(t) & d(t)
            \end{pmatrix}, \qquad t\geq0.
        \end{align}
        Observe that 
        \begin{align}
            \label{zeta1}
            & \sin(\zeta t)=\sin(\zeta_1t)\cosh(\zeta_2t)+i\cos(\zeta_1 t)\sinh(\zeta_2t),\\
            \label{zeta2}
            &\cos(\zeta t)=\cos(\zeta_1 t)\cosh(\zeta_2t)-i\sin(\zeta_1t)\sinh(\zeta_2t),
        \end{align} 
        where $\zeta=\zeta_1+i\zeta_2$, and $\zeta_2>0$.
        Then, $b(t)=-\sin(\zeta t)\neq0$ for $t>0$.
        We can compute
        \begin{equation}
            \frac{d(t)}{b(t)}=\frac{a(t)}{b(t)}=-\cot(\zeta t)=-\frac{\sin(2\zeta_1t)}{\cosh(2\zeta_2t)-\cos(2\zeta_1t)}+i\frac{\sinh(2\zeta_2t)}{\cosh(2\zeta_2t)-\cos(2\zeta_1t)}.
        \end{equation}
        Observe that $\Im(a(t)/b(t))>0$ for $t>0$. We can also use Theorem \ref{intro.main-2} to infer the integral representation of the propagator. This is, up to a phase, formula \eqref{propEx3}.

Moreover, the formulas above give
\[
    |b(t)|^2\,\Im\left(\frac{a(t)}{b(t)}\right)
    =\frac12\sinh(2\zeta_2t).
\]
Thus, by Theorem~\ref{thm-main3}~$(ii)$, the critical value
of $\gamma$ at time $t_0$ is $2/\sinh(2\zeta_2t_0)$.
For $\gamma>2/\sinh(2\zeta_2t_0)$ we get $u_0=0$. For $\gamma=2/\sinh(2\zeta_2t_0)$, a nonzero initial
$u_0\in L^p(\bR)$ is possible only if $p=\infty$ and $\Re b(t_0)\neq0$.
Since
\[
    \Re a(t_0)=\cos(\zeta_1t_0)\cosh(\zeta_2t_0),
    \qquad
    \Re b(t_0)=-\sin(\zeta_1t_0)\cosh(\zeta_2t_0),
\]
the latter condition is equivalent to $\sin(\zeta_1t_0)\neq0$.
In that case,
\[
    -\frac{\Re a(t_0)}{\Re b(t_0)}=\cot(\zeta_1t_0),
\]
and the asserted formula for $u_0$ follows from
\eqref{intro.def-fb20}.  
    \end{proof}
    \section{Concluding remarks}\label{sec:conclusions}
    The results above describe how the geometry of a positive symplectic matrix determines the critical decay threshold and the existence of extremizers in dimension one.
    A natural continuation concerns higher-dimensional oscillator semigroups. The unitary theory already shows that anisotropic Gaussian weights and degenerate symplectic blocks lead to directional uncertainty principles \cite{CGM2025,CGP2026}. Extending this description to positive complex symplectic matrices requires an analysis of the interaction between Gaussian damping and the relevant symplectic subspaces.
    
    A further direction concerns propagators generated by quadratic Hamiltonians with lower-order perturbations. The results in \cite{EKPV2008,EKPV2010,CassanoFanelli1} provide a natural point of comparison. In particular, bounded complex perturbations may alter Gaussian decay profiles \cite{CassanoFanelli1}.
    It is therefore natural to investigate which assumptions on the perturbation preserve uniqueness thresholds and which features of the critical regime survive.

\bibliographystyle{abbrv}
\bibliography{Bibliography}
\end{document}